\documentclass[11pt,a4paper]{article}

\usepackage{amsmath,amssymb}
\newcommand{\D}{\Delta}
\newcommand{\NN}{\frac{N}{p}}

\newcommand{\e}{\epsilon}
\newcommand{\va}{\varphi}
\newcommand{\n}{\nabla}
\newcommand{\N}{\frac{N}{2}}

\newcommand{\p}{\partial}

\newcommand{\R}{\mathbb{R}}
\newcommand{\T}{\mathbb{T}}

\newtheorem{definition}{Definition}
\newtheorem{theorem}{Theorem}

\newtheorem{proposition}{Proposition}
\newtheorem{corollaire}{Corollary}
\newtheorem{remarka}{Remark}

\newtheorem{lemme}{Lemma}

\title{Regularity of weak solutions of the compressible isentropic Navier-Stokes equation}
\author{Boris Haspot \thanks{Karls Ruprecht Universit\"at Heidelberg, Institut for Applied Mathematics,
Im Neuenheimer Feld 294,
D-69120 Heildelberg, Germany.
Tel. 49(0)6221-54-6112,
haspot@univ-paris12.fr
}}
\date{}
\begin{document}
\maketitle
\begin{abstract}
Regularity and uniqueness of weak solution of the compressible isentropic Navier-Stokes
equations is proven for small time in dimension $N=2,3$ under periodic boundary conditions.
In this paper, the initial density is not required to have a positive lower bound and the pressure
law is assumed to satisfy a condition that reduces to $\gamma>1$ when $N=2,3$ and $P(\rho)=a\rho^{\gamma}$.
In a second part we prove a condition of blow-up in slightly subcritical initial data when $\rho\in L^{\infty}$.
We finish by proving that weak solutions in $\T^{N}$ turn out to be smooth as long as the density remains
bounded in $L^{\infty}(L^{(N+1+\e)\gamma})$ with $\e>0$ arbitrary small.
\end{abstract}
\section{Introduction}
The Navier-Stokes equations are the basic model describing the evolution of a viscous compressible gas, As emphasized in many papers related to compressible fluid dynamics \cite{10,11,16,19,20}, vacuum is a major difficulty when trying to prove global existence and strong regularity results. As a matter of fact, starting from initial densities that have positive lower bounds, local existence of smooth solutions can be proved by classical means, since lower bounds on the density persists for small enough time. This paper is devoted to the proof of well-posedness results. We want to prove next that the norm $L^{p}(L^{q})$ on the pressure $P(\rho)$ control the breakdown of strong solutions of the Navier-Stokes equations when $N=2,3$. In other words, if a solution of of the Navier-Stokes equations is initially suitably smooth and loses it regularity at some later time, then the maximum norm of the density grows without bounds at the critical time approaches. Let us first recall the periodic compressible isentropic Navier-Stokes equations in $\T^{N}$ ($N\geq2$).
\begin{equation}
\begin{cases}
\begin{aligned}
&\p_{t}\rho+{\rm div}(\rho u)=0\;\;\;\mbox{in}\;\;{\cal D}^{'}((0,T)\times\T^{N}),\\
&\p_{t}(\rho u)+{\rm div}(\rho u\otimes u)-{\rm
div}(2\mu(\rho)Du)-\n(\lambda(\rho){\rm div}u)
+\n P(\rho)=\rho g,
\end{aligned}
\end{cases}
\label{1}
\end{equation}
The unknowns $\rho$, $u$ respectively correspond to the density of the gas $\rho\geq0$ and its velocity $u\in\R^{N}$. The last equation of (\ref{1}) defines the pressure $P$ which is assumed to be increasing in $\rho$ for physical reasons. Usually $P(\rho)=P_{\gamma,a}=a\rho^{\gamma}$ for some positive constant $a$ and some $\gamma\geq1$. The viscosity coefficients ares assumed to satisfy $\mu>0$, $N\lambda+2\mu\geq0$, and the external forces $g$ to belong to $L^{2}((0,T)\times\T^{N})^{N}$ for all $T>0$. Finally we complement the above system with initial conditions
\begin{equation}
\begin{cases}
\begin{aligned}
&\rho_{\ t=0}=\rho_{0}\geq 0,\\
&\rho u_{\ t=0}=m_{0}.
\end{aligned}
\end{cases}
\label{2}
\end{equation}
Before the remarkable work of P-L Lions, very little was known about solutions of the compressible isentropic Navier-Stokes equation at least when $N\geq 2$. In \cite{13} he proved a global existence theorem and weak stability results for $P_{\gamma,a}$ pressure laws under the following assumptions on the initial data
\begin{equation}
\begin{cases}
\begin{aligned}
&\rho_{0}\in L^{1}(\T^{N})\cap L^{\gamma}(\T^{N}), \rho_{0}\geq0,\\
&\frac{m_{0}^{2}}{\rho_{0}}\in L^{1}(\T^{N}),
\end{aligned}
\end{cases}
\label{3}
\end{equation}
where we agree that $\frac{m_{0}^{2}}{\rho_{0}}=0$ on $\{x\in\T^{N}\;\mbox{such that}\;\rho_{0}(x)=0\}$. More precisely he proved
\begin{theorem}
\label{theorem1}
We assume (\ref{3}) and $\gamma>1$ if $N=2$, $\gamma>\frac{3}{2}$ if $N=3$. Then there exists a solution $(\rho,u)\in L^{\infty}(0,\infty;L^{\gamma}(\T^{N}))\times L^{2}(0,\infty;H^{1}(\T^{N}))^{N}$ satisfying in addition $\rho\in C([0,\infty),L^{p}(\T^{N}))$ if $1\leq p<\gamma$, $\rho|u|^{2}\in L^{\infty}(0,\infty;L^{1}(\T^{N}))$, $\rho\in L^{q}_{loc}([0,\infty);L^{q}(\T^{N}))$ for $1\leq q\leq \gamma-1+\frac{2}{N}\gamma$. Moreover, when $f=0$, for almost all $t\geq 0$, we have
\begin{equation}
\begin{aligned}
&\int_{\T^{N}}(\frac{1}{2}\rho|u|^{2}+\frac{a}{\gamma-1}\rho^{\gamma})(t,x)dx+\int^{t}_{0}ds
\int_{\T^{N}}(\mu|\n u|^{2}+(\lambda+\mu)({\rm div}u)^{2})dxds\\
&\hspace{7cm}\leq \int_{\T^{N}}(\frac{1}{2}\rho_{0}|u_{0}|^{2}+\frac{a}{\gamma-1}\rho_{0}^{\gamma})(x)dx.
\end{aligned}
\label{4}
\end{equation}
\end{theorem}
In addition, he proved similar results for more general pressure laws $P(\rho)$ such that
\begin{equation}
\begin{cases}
\begin{aligned}
&\int^{1}_{0}\frac{P(s)}{s^{2}}ds<+\infty,\\
&\lim\inf\rightarrow_{s\rightarrow+\infty}\frac{P(s)}{s^{\gamma}}>0,
\end{aligned}
\end{cases}
\label{e4}
\end{equation}
for some $\gamma$ satisfying the above condition of theorem \ref{theorem1}. Notice that the main difficulty for proving Lions' theorem consists
in exhibiting strong compactness properties of the density $\rho$ in
$L^{p}_{loc}(\R^{+}\times\R^{N})$ spaces required to pass to the
limit in the pressure
term $P(\rho)=a\rho^{\gamma}$. Let us mention that Feireisl and his collaborators in \cite{5F1,5F2,5F3} generalized the result to
any $\gamma>\frac{N}{2}$ for $N\geq 2$ in establishing that we can obtain renormalized
solution without imposing that $\rho\in
L^{2}_{loc}(\R^{+}\times\R^{N})$ (a property that was  needed in Lions' approach in dimension
$N=2,3$ giving the further condition  $\gamma-1+\frac{2\gamma}{N}\geq2$), for this he
introduces the concept of oscillation defect measure evaluating the
loss of compactness.\\
In \cite{5BD1} Bresch and Desjardins show a result of global
existence of weak solution for the non isothermal Navier-Stokes
system assuming density dependence of $\mu(\rho)$ and $\lambda(\rho)$, considering perfect gas law with some cold pressure close to the vacuum, and the following relation:
\begin{equation}
\lambda(\rho)=2(\rho\mu^{'}(\rho)-\mu(\rho)). \label{5coeff}
\end{equation}
The key point in this paper is to
show that the structure of the diffusion term provides some regularity for the
density thanks to a new mathematical entropy inequality. This one has been discovered in \cite{5BD2}, we call it the BD entropy.
Mellet and Vasseur by using the BD entropy, get in \cite{5MV} a very interesting stability result.
The interest of this result is to
consider conditions where the viscosity coefficients vanish on the
vacuum set. It includes the case $\mu(\rho)=\rho$, $\lambda(\rho)=0$
(when $N=2$ and $\gamma=2$, where we recover the Saint-Venant model
for Shallow water). The key to the proof is  a new energy inequality
on the velocity and a gain of integrability, which
allows to pass to the limit. Unfortunately, the construction of approximate solutions satisfying: energy estimates, BD mathematical entropy and Mellet-Vasseur estimates is far from being proven except in dimension one or with symmetry assumptions, see \cite{aMV1}, \cite{aJJX}, \cite{aGJX}. Note that approximate solutions construction process has been proposed in \cite{aBD} satisfying energy estimates and BD mathematical entropy. This means that only global existence of weak solutions with some extra terms or cold pressure exists in dimension greater than 2.\\
The existence and uniqueness of local classical solutions for (\ref{1})
with smooth initial data such that the density $\rho_{0}$ is bounded
and bounded away from zero (i.e.,
$0<\underline{\rho}\leq\rho_{0}\leq M$)
has been stated by Nash in \cite{5Na}. Let us emphasize that no stability condition was required there.\\
On the other hand, for small smooth perturbations of a stable
equilibrium with constant positive density, global well-posedness
has been proved in \cite{5MN}. Many works in the case of the one dimension have been devoted
to the qualitative behavior of solutions for large time (see for
example \cite{9,11}). Refined functional analysis has been used
for the last decades, ranging from Sobolev, Besov, Lorentz and
Triebel spaces to describe the regularity and long time behavior of
solutions to the compressible model \cite{5So}, \cite{20},
\cite{5H4}, \cite{5HZ}, \cite{5K1}. \\
The use of critical functional frameworks led to several new well-posedness results for compressible
fluids. In addition to have a norm invariant by (\ref{1}),
appropriate functional space for solving (\ref{1}) must provide a control on the $L^{\infty}$
norm of the density (in order to avoid vacuum and loss of ellipticity). For that reason,
we restricted our study to the case where the initial data $(\rho_{0},u_{0})$ and external force $f$
are such that, for some positive constant $\bar{\rho}$:
$$(\rho_{0}-\bar{\rho})\in B^{\NN}_{p,1},\;u_{0}\in B^{\frac{N}{p_{1}}-1}_{p_{1},1}\;\;\mbox{and}\;\;f\in L^{1}_{loc}(\R^{+},\in B^{\frac{N}{p_{1}}-1}_{p_{1},1})$$
with $(p,p_{1})\in [1,+\infty[$ good choose.
The most important result come from R. Danchin in
\cite{DL} which show the existence of global solution
and uniqueness with initial data close from a equilibrium, and he
obtains a similar result in finite time. The interest is that he works
in {\it critical} Besov space ({\it critical} in the sense
of the scaling of the equation). More precisely to speak roughly, he get strong solution with
initial data in $B^{\frac{N}{2}}_{2,1}\cap B^{\frac{N}{2}-1}_{2,1}\times (B^{\frac{N}{2}-1}_{2,1})^{N}$. Here compared with
the result on Navier-Stokes incompressible, he needs to control the vacuum and the norm $L^{\infty}$ of the density in the goal
to use the parabolicity of the momentum equation and to have some properties of multiplier spaces.
That's why Danchin works in Besov spaces with a third index $r=1$ for the density, and it's the same for the velocity as the equations are linked. In \cite{DW}, R. Danchin generalize the previous result with large initial data on the density.\\
In \cite{DW}, however, we hand to have $p=p_{1}$, indeed in this article there exists a very strong coupling between the pressure
and the velocity. To be more precise, the pressure term is considered as a term of rest for the elliptic operator in the momentum equation of (\ref{1}). This paper improve the results of R. Danchin in \cite{DL,DW}, in the sense that the initial density belongs to larger spaces $B^{\NN}_{p,1}$ with $p\in[1,+\infty[$.
The main idea of this paper is to introduce a new variable than the velocity in the goal to \textit{kill}
the relation of coupling between the velocity and the density. In \cite{CD}, F. Charve and R. Danchin and in \cite{CMZ} Q. Chen et al
generalize the results from \cite{DG} by choosing more general initial data. In particular they works with general Besov space constructed on $L^{p}$, however they added some conditions on $p$ ($p<2N$) to get global solutions. For results of strong solutions with general viscosity coefficients we refer to \cite{CMZ1,H,H1}.\\
In \cite{H3}, we address the question of local well-posedness in the critical functional framework under the assumption
that the initial density belongs to critical Besov space with a index of integrability different of this of the velocity.
We adapt the spirit of the results of \cite{AP} and \cite{H} which treat the case of Navier-Stokes incompressible with dependent density (at the difference than in these works the velocity and the density are naturally decoupled).
The main idea of this paper is to introduce a new variable than the velocity in the goal to ''\textit{kill}''
the coupling between the velocity and the density.
We introduce a new variable $v_{1}$ to control the velocity where to avoid the coupling between the density and the velocity, we analyze by a new way the pressure term (in particular we will use this variable $v_{1}$ in this article). This idea is inspired from the works of D. Hoff, P-L Lions and D. Serre about the the famous effective pressure.
 More precisely we write the gradient of the pressure as a Laplacian of the variable $v_{1}$, and we introduce this term in the linear part of the momentum equation. We have then a control on $v_{1}$ which can write roughly as $u-{\cal G}P(\rho)$ where ${\cal G}$ is a pseudodifferential operator of order $-1$. By this way, we have canceled the coupling between $v_{1}$ and the density, we next verify easily that we have a control Lipschitz of the gradient of $u$ (it is crucial to estimate the density by the transport equation).
This result allows us to reach some critical initial data in the sense that we are not very far to choose $(\rho_{0}-\bar{\rho},u_{0})$ in $B^{0}_{\infty,1}\times B^{0}_{N,1}$.
\\
\\
On the other hand there have been few existence results on the strong solutions for the general case of nonnegative initial densities. The first result was proved by R. Salvi and I. Straskraba. They showed in \cite{CCK17} that
if $\Omega$ is a bounded domain, $P=P(\cdot)\in C^{2}[0,\infty)$, $\rho_{0}\in H^{2}$, $u_{0}\in H^{1}_{0}\cap H^{2}$ and the compatibility condition:
\begin{equation}
Lu_{0}+\n P(\rho_{0})=\rho_{0}^{\frac{1}{2}}g,\;\;\mbox{for some}\;g\in L^{2},
\label{1.3}
\end{equation}
is satisfied, then there exists a unique local strong solution $(\rho,u)$ to the initial boundary value problem (\ref{1}). H. J. Choe and H. Kim proved in \cite{CCK3} a similar existence result when $\Omega$ is either a bounded domain or the whole space, $P(\rho)=a\rho^{\gamma}$ ($a>0$, $\gamma>1$), $\rho_{0}\in L^{1}\cap H^{1}\cap W^{1,6}$, $u_{0}\in D^{1}_{0}\cap D^{2}$ and the condition (\ref{1.3}) is satisfied.\\
B. Desjardins in \cite{CCK5} proved the local existence of a weak solution solution $(\rho,u)$ with a bounded nonnegative density to the periodic boundary value problem (\ref{1}) as long as $\sup_{0\leq t\leq T^{*}}(\|\rho(t)\|_{L^{\infty}(\T^{3})}+\|\n u(t)\|_{L^{2}(\T^{3})})<+\infty$.\\
This paper is devoted to improve the works in \cite{CCK5} and \cite{CCK3} by choosing very low regularity on the velocity. In the sequel we will note $\frac{d}{dt}=\p_{t}+u\cdot\n$ and $\dot{f}=\frac{d}{dt}f$.\\
The viscosity coefficients will suppose constant in the sequel and are assumed to satisfy:
\begin{equation}
\mu>0,\;\;0<\lambda\frac{5}{4}\mu.
\label{condviscosity}
\end{equation}
It follow that there is a $l>4$, which will be fixed throughout, such that:
\begin{equation}
\frac{\mu}{\lambda}>\frac{(l-2)^{2}}{4(l-1)}.
\label{condviscosity1}
\end{equation}
In the sequel we will assume that $g\in E^{1}_{T}$ with:
$$\|g\|_{E^{1}_{T}}=\|g\|_{L^{\infty}_{T}(L^{2})}+\|g\|_{L^{2}_{T}(L^{2})}+\int^{T}_{0}(f(s)^{7}\|\n f\|_{L^{4}}^{2}+\int^{T}_{0}\int_{\R^{N}}f(s)^{5}|\p_{s}f|^{2}dsdx,$$
where $f(s)=\min(1,s)$ .
We obtain now the following existence of weak solutions in finite:
\begin{theorem}
\label{theo1}
Let $N=2,3$. Assume that $\mu$ and $\lambda$ verify (\ref{condviscosity}) and (\ref{condviscosity1}). We assume that $\rho_{0}\in L^{\infty}(\R^{N})$, $\rho_{0}^{\frac{1}{4}}u_{0}\in L^{4}$ and $\rho_{0}^{\frac{1}{2}}u_{0}\in L^{2}$.
Moreover $g$ is in $E^{1}_{T}$ and $g\in L^{4}(L^{4})\cap L^{\infty}(L^{4})$.
\begin{itemize}
 \item  There exists $T_{0}\in (0,+\infty]$ and a weak solution $(\rho,u)$ to the system (\ref{1}) in $[0,T_{0}]$ such that for all $T<T_{0}$ and with $f(t)=\min(t,1)$:
\begin{equation}
 \begin{aligned}
&\sup_{0<t\leq T_{0}}\int_{\R^{N}}[\frac{1}{2}\rho(t,x)|u(t,x)|^{2}+|P(\rho(t,x))|+f(t)|\n u(t,x)|^{2}dx\\
&+\sup_{0<t\leq T_{0}}\int_{\R^{N}}[\frac{1}{2}\rho(t,x)f(t)^{N}(\rho|\dot{u}(t,x)|^{2}+|\n \omega(t,x)|^{2})dx\\
&+\int^{+\infty}_{0}\int_{\R^{N}}[|\n u|^{2}+f(s)\rho|\dot{u}|^{2}+|\omega|^{2})+\sigma^{N}|\n \dot{u}|^{2}]dxdt\leq C C_{0}.
 \end{aligned}
\label{1.21}
\end{equation}
\begin{equation}
\begin{aligned}
&\sup_{0<t\leq T_{0}}\int_{\R^{N}}\rho(t,x)|u(t,x)|^{4}dx\leq C_{0}.
\end{aligned}
 \label{1.211}
\end{equation}

where $C_{0}$ depends of the initial data $\rho_{0}$ and $u_{0}$. We obtain moreover:
\begin{equation}
\rho\in L^{\infty}((0,T)\times\R^{N})\cap L^{\infty}_{T}(L^{\gamma}_{2}(\R^{N})),
 \label{12}
\end{equation}
\begin{equation}
\begin{cases}
 \begin{aligned}
&\sqrt{\rho}\p_{t}u\in L^{2}_{t}(L^{2}(\R^{N})),\\
&\sqrt{t}{\cal P}u\in L^{2}_{T}(H^{2}(\R^{N})),\\
&\sqrt{t}G=\sqrt{t}[(\lambda+2\mu){\rm div}u-P(\rho)]\in L^{2}_{T}(H^{1}(\R^{N})),\\
&\sqrt{t}\n u\in L^{\infty}_{T}(L^{2}(\R^{N})),
 \end{aligned}
\end{cases}
\label{13}
\end{equation}
where ${\cal P}$ denotes the projection on the space of divergence-free vector fields, and
\item In addition if we assume that $u_{0}\in H^{\frac{N}{2}-1+\e}$ with $\e>0$, we obtain the following estimates:
\begin{equation}
 \begin{aligned}
&\sup_{t>0}\int_{\R^{N}}[\frac{1}{2}\rho(t,x)|u(t,x)|^{2}+|P(\rho(t,x))|+t^{2-\N-\e}|\n u(t,x)|^{2}dx\\
&+\sup_{0<t\leq T_{0}}\int_{\R^{N}}[\frac{1}{2}\rho(t,x)t^{\sigma}(\rho|\dot{u}(t,x)|^{2}+|\n \omega(t,x)|^{2})dx\\
&+\int^{+\infty}_{0}\int_{\R^{N}}[|\n u|^{2}+t^{2-\N-\e}|\dot{u}|^{2}+t^{\sigma}|\n \dot{u}|^{2}]dxdt\leq C(C_{0}+C_{f})^{\theta},
 \end{aligned}
\label{b1.22}
\end{equation}
where:
\begin{equation}
 \begin{cases}
  \begin{aligned}
&\sigma=3-\N-\e,\;\;\;\mbox{if}\;\;N=2,\\
&\sigma=\max(3-\N-\e,6-\frac{3}{2}N-3\e),\;\;\;\mbox{if}\;\;N=3,
  \end{aligned}
\end{cases}
\label{c24}
\end{equation}
and:
\begin{equation}
\sup_{0\leq t\leq T_{0}}\|u(t,\cdot)\|_{H^{\N-1+\e}}\leq CC_{0}^{\theta}\;\;\mbox{and}\;\;
\sup_{0\leq t\leq T_{0}}\|\rho^{\frac{1}{r}}u(t,\cdot)\|_{L^{r}}\leq C(r)C_{0}^{\theta},
 \label{b1.23}
\end{equation}
with $r\in(1,4]$ and:
\begin{equation}
\n u\in L^{1}_{T}(BMO).
\label{c1.23}
\end{equation}
\item The regularity properties (\ref{1.21}), (\ref{13}), (\ref{b1.22}), (\ref{b1.23}) and (\ref{c1.23}) hold as long as:
\begin{equation}
\sup_{t\in[0,T]}\|\rho\|_{L^{\infty}_{t}(L^{\infty}(\R^{N}))}<+\infty.
\label{11.23}
\end{equation}
\end{itemize}
\end{theorem}
\begin{remarka}
\begin{itemize}
 \item In this theorem we do not assume hypothesis on the vacuum as in \cite{H1} where $\frac{1}{\rho_{0}}\in L^{\infty}$. However we get a control $L^{1}(BMO)$ on the gradient of the velocity. We know that for incompressible Navier-Stokes this hypothesis is enough to get uniqueness. In this sense we can consider our result as a theorem of strong solutions in a weak sense.
\item In this result we improve the works of B. Desjardins in \cite{CCK5} by the fact that we can choose more general initial data. The second important point is that estimates (\ref{1.21}) and (\ref{1.211}) can hold as long as $\rho\in L^{\infty}$ in dimension $2$ and $3$.
\end{itemize}
 \end{remarka}
Let $\bar{\rho}>0$. In the sequel we will note $q=\rho-\bar{\rho}$.
In the following corollary we obtain strong solution by adding some regularity on the initial density and we assume that $\rho_{0}$ is away from the vacuum. The goal is to get by this extra regularity on the density a control Lipschitz of $\n u$.
\begin{corollaire}
\label{corollaire1}
Under the hypothesis of theorem \ref{theo1}, we assume in addition that there exists $c>0$ such that $\rho_{0}\geq c$. Moreover we assume that $q_{0}\in B^{\e}_{\infty,\infty}$ if $P(\rho)=K\rho$ and $q_{0}\in B^{1+\e}_{N,\infty}$ if $P$ is a general pressure.
The solutions of theorem \ref{theo1} are then unique and verify locally in time (\ref{1.21}), (\ref{13}), (\ref{b1.22}), (\ref{b1.23}) and:
$$\n u\in L^{1}_{T_{0}}(L^{\infty}).$$
Moreover if $\rho$ in $L^{\infty}(\R\times\T^{N})$ then the solution is global.
\end{corollaire}
In the following theorem we want improve the criterion of blow-up of corollary \ref{corollaire1}. More precisely we prove that it is just necessary to control the norm $L^{\infty}(L^{(N+1+\e)\gamma})$ with $\e>0$ of the density when $P(\rho)=a\rho^{\gamma}$ with $\gamma\geq1$ to obtain global strong solution. This result is to connect with the works of Serrin for incompressible Navier-Stokes system where in the compressible Navier-Stokes, the pressure plays the role of the velocity.
\begin{theorem}
\label{theo3}
Let $\lambda=0$ and $g$ as in theorem \ref{theo1} and $g\in L^{\infty}(L^{\infty})$. Let $P(\rho)=a\rho^{\gamma}$ with $a>0$ and $\gamma\geq 1$. Assume that $(\rho_{0},u_{0})\in (L^{\gamma}_{2}\cap L^{\infty}\cap B^{1+\e}_{N,\infty})\times (L^{2}\cap L^{\infty}\cap H^{\N-1+\e})$ with $\e>0$. Moreover $\rho_{0}$ check $\rho_{0}\geq c>0$.\\
Let $(\rho,u)$ a weak solution of system (\ref{1}) on $[0,T)$ with the previous initial data which satisfies the following condition:
\begin{equation}
\rho\in L^{\infty}(0,+\infty,L^{(N+1+\e)\gamma}),
\label{egcrucial}
\end{equation}
with $\e>0$. Then $(\rho,u)$ is unique and verify locally in time  (\ref{1.21}), (\ref{13}), (\ref{b1.22}), (\ref{b1.23}) and:
$$\n u\in L^{1}_{loc}(L^{\infty}).$$
\end{theorem}
\begin{remarka}
\begin{itemize}
\item In this theorem we can see that by compare with the incompressible Navier-Stokes equation, the good variable to control is not the velocity but the pressure. Indeed if we control enough the pressure, we get integrability on the velocity.
 \item In this theorem we need to assume hypothesis on the vacuum as in \cite{H1} where $\frac{1}{\rho_{0}}\in L^{\infty}$ in the goal to get a control Lipschitz on the velocity.
\item This initial data are considered slightly surcritical in dimension $3$ except that $u_{0}\in L^{\infty}$ (it is probably possible to improve this fact).
\item Here we need to assume that $\lambda=0$ to get a control $L^{\infty}$ on $u$ as in the article of A. Mellet and A. Vasseur in \cite{5MV2}. We recall that in this paper they need of a control on $P(\rho)\in L^{\infty}(L^{3+\e})$ with $\e>0$ for $N=3$ by using some De Giorgi technics used by A. Vasseur in \cite{V6} to reprove the famous result of Caffarelli-Kohn and Nirenberg in \cite{CKN}. As in \cite{5MV2}, the pressure plays a important role, and in some sense the pressure is the good variable to control to get global strong solution. The pressure plays the role of the velocity for Navier-Stokes incompressible when we have compressible Navier-Stokes system where a structure of type effective pressure exists.
\item We can justify the previous statement by seeing the results of D. Bresch and B. Desjardins in \cite{5BD}. They have a new mathematical entropy for a class of such viscosity coefficients which gives some norms on the gradient of $\rho$ (see \cite{5BD}), in particular we can obtain the relation (\ref{egcrucial}). However this type of viscosity coefficient kill the structure of effective pressure and we can apply our proof. It's again a proof that the structure of the viscosity coefficients plays a crucial role for compressible Navier-Stokes system.
\item As in the case of incompressible Navier-Stokes system (see \cite{CKN}), there is a big gap between obtaining (\ref{egcrucial}) and by this way have a control on $\rho$ in $L^{\infty}$.
\item V. A. Waigant has builded in \cite{CCK25} explicit solutions for which the maximal integrability of the density correspond to $L^{q}(0,1,L^{q})$ with $q=\frac{\gamma(3N+2)-N}{2N}$. It means that
(\ref{egcrucial}) fails in this case except that the force term is less regular than in our case. It means that the regularity of $g$ is crucial to get strong solution.
\end{itemize}
 \end{remarka}
\begin{remarka}
We believe that our method can be adapted to the euclidian space $\mathbb{R}^{N}$. This is the object of our future work.
\end{remarka}
Our paper is structured as follows. In section \ref{section2}, we give a few notation and briefly introduce the basic Fourier analysis
techniques needed to prove our result. In section \ref{section3} and \ref{section4}, we prove  a priori estimate on the density and the velocity. In section \ref{section5} and section \ref{section6}, we prove the theorem \ref{theo1} and corollary \ref{corollaire1}. We finish in the section \ref{section7} by the proof of theorem \ref{theo3}.
\section{Littlewood-Paley theory and Besov spaces}
\label{section2}
Throughout the paper, $C$ stands for a constant whose exact meaning depends on the context. The notation $A\lesssim B$ means
that $A\leq CB$.
For all Banach space $X$, we denote by $C([0,T],X)$ the set of continuous functions on $[0,T]$ with values in $X$.
For $p\in[1,+\infty]$, the notation $L^{p}(0,T,X)$ or $L^{p}_{T}(X)$ stands for the set of measurable functions on $(0,T)$
with values in $X$ such that $t\rightarrow\|f(t)\|_{X}$ belongs to $L^{p}(0,T)$.
Littlewood-Paley decomposition  corresponds to a dyadic
decomposition  of the space in Fourier variables.
We can use for instance any $\varphi\in C^{\infty}(\T^{N})$,
supported in
${\cal{C}}=\{\xi\in\T^{N}/\frac{3}{4}\leq|\xi|\leq\frac{8}{3}\}$
such that:
$$\sum_{l\in\mathbb{Z}}\varphi(2^{-l}\xi)=1\,\,\,\,\mbox{if}\,\,\,\,\xi\ne 0.$$
Denoting $h={\cal{F}}^{-1}\varphi$, we then define the dyadic
blocks by:
$$\D_{l}u=\varphi(2^{-l}D)u=2^{lN}\int_{\T^{N}}h(2^{l}y)u(x-y)dy\,\,\,\,\mbox{and}\,\,\,S_{l}u=\sum_{k\leq
l-1}\D_{k}u\,.$$ Formally, one can write that:
$$u=\sum_{k\in\mathbb{Z}}\D_{k}u\,.$$
This decomposition is called homogeneous Littlewood-Paley
decomposition. Let us observe that the above formal equality does
not hold in ${\cal{S}}^{'}(\T^{N})$ for two reasons:
\begin{enumerate}
\item The right hand-side does not necessarily converge in
${\cal{S}}^{'}(\T^{N})$.
\item Even if it does, the equality is not
always true in ${\cal{S}}^{'}(\T^{N})$ (consider the case of the
polynomials).
\end{enumerate}
\subsection{Homogeneous Besov spaces and first properties}
\begin{definition}
For
$s\in\R,\,\,p\in[1,+\infty],\,\,q\in[1,+\infty],\,\,\mbox{and}\,\,u\in{\cal{S}}^{'}(\T^{N})$
we set:
$$\|u\|_{B^{s}_{p,q}}=(\sum_{l\in\mathbb{Z}}(2^{ls}\|\D_{l}u\|_{L^{p}})^{q})^{\frac{1}{q}}.$$
The Besov space $B^{s}_{p,q}$ is the set of temperate distribution $u$ such that $\|u\|_{B^{s}_{p,q}}<+\infty$.
\end{definition}
\begin{remarka}The above definition is a natural generalization of the
nonhomogeneous Sobolev and H$\ddot{\mbox{o}}$lder spaces: one can show
that $B^{s}_{\infty,\infty}$ is the nonhomogeneous
H$\ddot{\mbox{o}}$lder space $C^{s}$ and that $B^{s}_{2,2}$ is
the nonhomogeneous space $H^{s}$.
\end{remarka}
\begin{proposition}
\label{derivation,interpolation}
The following properties holds:
\begin{enumerate}
\item there exists a constant universal $C$
such that:\\
$C^{-1}\|u\|_{B^{s}_{p,r}}\leq\|\n u\|_{B^{s-1}_{p,r}}\leq
C\|u\|_{B^{s}_{p,r}}.$
\item If
$p_{1}<p_{2}$ and $r_{1}\leq r_{2}$ then $B^{s}_{p_{1},r_{1}}\hookrightarrow
B^{s-N(1/p_{1}-1/p_{2})}_{p_{2},r_{2}}$.
\item $B^{s^{'}}_{p,r_{1}}\hookrightarrow B^{s}_{p,r}$ if $s^{'}> s$ or if $s=s^{'}$ and $r_{1}\leq r$.
\end{enumerate}
\label{interpolation}
\end{proposition}
Before going further into the paraproduct for Besov spaces, let us state an important proposition.
\begin{proposition}
Let $s\in\R$ and $1\leq p,r\leq+\infty$. Let $(u_{q})_{q\geq-1}$ be a sequence of functions such that
$$(\sum_{q\geq-1}2^{qsr}\|u_{q}\|_{L^{p}}^{r})^{\frac{1}{r}}<+\infty.$$
If $\mbox{supp}\hat{u}_{1}\subset {\cal C}(0,2^{q}R_{1},2^{q}R_{2})$ for some $0<R_{1}<R_{2}$ then $u=\sum_{q\geq-1}u_{q}$ belongs to $B^{s}_{p,r}$ and there exists a universal constant $C$ such that:
$$\|u\|_{B^{s}_{p,r}}\leq C^{1+|s|}\big(\sum_{q\geq-1}(2^{qs}\|u_{q}\|_{L^{p}})^{r}\big)^{\frac{1}{r}}.$$
\label{resteimp1}
\end{proposition}
Let now recall a few product laws in Besov spaces coming directly from the paradifferential calculus of J-M. Bony
(see \cite{5BJM}) and rewrite on a generalized form in \cite{AP} by H. Abidi and M. Paicu (in this article the results are written
in the case of homogeneous spaces but it can easily generalize for the nonhomogeneous Besov spaces).
\begin{proposition}
\label{produit1}
We have the following laws of product:
\begin{itemize}
\item For all $s\in\R$, $(p,r)\in[1,+\infty]^{2}$ we have:
\begin{equation}
\|uv\|_{B^{s}_{p,r}}\leq
C(\|u\|_{L^{\infty}}\|v\|_{B^{s}_{p,r}}+\|v\|_{L^{\infty}}\|u\|_{B^{s}_{p,r}})\,.
\label{2.2}
\end{equation}
\item Let $(p,p_{1},p_{2},r,\lambda_{1},\lambda_{2})\in[1,+\infty]^{2}$ such that:$\frac{1}{p}\leq\frac{1}{p_{1}}+\frac{1}{p_{2}}$,
$p_{1}\leq\lambda_{2}$, $p_{2}\leq\lambda_{1}$, $\frac{1}{p}\leq\frac{1}{p_{1}}+\frac{1}{\lambda_{1}}$ and
$\frac{1}{p}\leq\frac{1}{p_{2}}+\frac{1}{\lambda_{2}}$. We have then the following inequalities:\\
if $s_{1}+s_{2}+N\inf(0,1-\frac{1}{p_{1}}-\frac{1}{p_{2}})>0$, $s_{1}+\frac{N}{\lambda_{2}}<\frac{N}{p_{1}}$ and
$s_{2}+\frac{N}{\lambda_{1}}<\frac{N}{p_{2}}$ then:
\begin{equation}
\|uv\|_{B^{s_{1}+s_{2}-N(\frac{1}{p_{1}}+\frac{1}{p_{2}}-\frac{1}{p})}_{p,r}}\lesssim\|u\|_{B^{s_{1}}_{p_{1},r}}
\|v\|_{B^{s_{2}}_{p_{2},\infty}},
\label{2.3}
\end{equation}
when $s_{1}+\frac{N}{\lambda_{2}}=\frac{N}{p_{1}}$ (resp $s_{2}+\frac{N}{\lambda_{1}}=\frac{N}{p_{2}}$) we replace
$\|u\|_{B^{s_{1}}_{p_{1},r}}\|v\|_{B^{s_{2}}_{p_{2},\infty}}$ (resp $\|v\|_{B^{s_{2}}_{p_{2},\infty}}$) by
$\|u\|_{B^{s_{1}}_{p_{1},1}}\|v\|_{B^{s_{2}}_{p_{2},r}}$ (resp $\|v\|_{B^{s_{2}}_{p_{2},\infty}\cap L^{\infty}}$),
if $s_{1}+\frac{N}{\lambda_{2}}=\frac{N}{p_{1}}$ and $s_{2}+\frac{N}{\lambda_{1}}=\frac{N}{p_{2}}$ we take $r=1$.
\\
If $s_{1}+s_{2}=0$, $s_{1}\in(\frac{N}{\lambda_{1}}-\frac{N}{p_{2}},\frac{N}{p_{1}}-\frac{N}{\lambda_{2}}]$ and
$\frac{1}{p_{1}}+\frac{1}{p_{2}}\leq 1$ then:
\begin{equation}
\|uv\|_{B^{-N(\frac{1}{p_{1}}+\frac{1}{p_{2}}-\frac{1}{p})}_{p,\infty}}\lesssim\|u\|_{B^{s_{1}}_{p_{1},1}}
\|v\|_{B^{s_{2}}_{p_{2},\infty}}.
\label{2.4}
\end{equation}
If $|s|<\NN$ for $p\geq2$ and $-\frac{N}{p^{'}}<s<\NN$ else, we have:
\begin{equation}
\|uv\|_{B^{s}_{p,r}}\leq C\|u\|_{B^{s}_{p,r}}\|v\|_{B^{\NN}_{p,\infty}\cap L^{\infty}}.
\label{2.5}
\end{equation}
\end{itemize}
\end{proposition}
\begin{remarka}
In the sequel $p$ will be either $p_{1}$ or $p_{2}$ and in this case $\frac{1}{\lambda}=\frac{1}{p_{1}}-\frac{1}{p_{2}}$
if $p_{1}\leq p_{2}$, resp $\frac{1}{\lambda}=\frac{1}{p_{2}}-\frac{1}{p_{1}}$
if $p_{2}\leq p_{1}$.
\end{remarka}
\begin{corollaire}
\label{produit2}
Let $r\in [1,+\infty]$, $1\leq p\leq p_{1}\leq +\infty$ and $s$ such that:
\begin{itemize}
\item $s\in(-\frac{N}{p_{1}},\frac{N}{p_{1}})$ if $\frac{1}{p}+\frac{1}{p_{1}}\leq 1$,
\item $s\in(-\frac{N}{p_{1}}+N(\frac{1}{p}+\frac{1}{p_{1}}-1),\frac{N}{p_{1}})$ if $\frac{1}{p}+\frac{1}{p_{1}}> 1$,
\end{itemize}
then we have if $u\in B^{s}_{p,r}$ and $v\in B^{\frac{N}{p_{1}}}_{p_{1},\infty}\cap L^{\infty}$:
$$\|uv\|_{B^{s}_{p,r}}\leq C\|u\|_{B^{s}_{p,r}}\|v\|_{B^{\frac{N}{p_{1}}}_{p_{1},\infty}\cap L^{\infty}}.$$
\end{corollaire}
The study of non stationary PDE's requires space of type $L^{\rho}(0,T,X)$ for appropriate Banach spaces $X$. In our case, we
expect $X$ to be a Besov space, so that it is natural to localize the equation through Littlewood-Payley decomposition. But, in doing so, we obtain
bounds in spaces which are not type $L^{\rho}(0,T,X)$ (except if $r=p$).
We are now going to
define the spaces of Chemin-Lerner in which we will work, which are
a refinement of the spaces
$L_{T}^{\rho}(B^{s}_{p,r})$.
$\hspace{15cm}$
\begin{definition}
Let $\rho\in[1,+\infty]$, $T\in[1,+\infty]$ and $s_{1}\in\R$. We set:
$$\|u\|_{\widetilde{L}^{\rho}_{T}(B^{s_{1}}_{p,r})}=
\big(\sum_{l\in\mathbb{Z}}2^{lrs_{1}}\|\D_{l}u(t)\|_{L^{\rho}(L^{p})}^{r}\big)^{\frac{1}{r}}\,.$$
We then define the space $\widetilde{L}^{\rho}_{T}(B^{s_{1}}_{p,r})$ as the set of temperate distribution $u$ over
$(0,T)\times\T^{N}$ such that 
$\|u\|_{\widetilde{L}^{\rho}_{T}(B^{s_{1}}_{p,r})}<+\infty$.
\end{definition}
We set $\widetilde{C}_{T}(\widetilde{B}^{s_{1}}_{p,r})=\widetilde{L}^{\infty}_{T}(\widetilde{B}^{s_{1}}_{p,r})\cap
{\cal C}([0,T],B^{s_{1}}_{p,r})$.
Let us emphasize that, according to Minkowski inequality, we have:
$$\|u\|_{\widetilde{L}^{\rho}_{T}(B^{s_{1}}_{p,r})}\leq\|u\|_{L^{\rho}_{T}(B^{s_{1}}_{p,r})}\;\;\mbox{if}\;\;r\geq\rho
,\;\;\;\|u\|_{\widetilde{L}^{\rho}_{T}(B^{s_{1}}_{p,r})}\geq\|u\|_{L^{\rho}_{T}(B^{s_{1}}_{p,r})}\;\;\mbox{if}\;\;r\leq\rho
.$$
\begin{remarka}
It is easy to generalize proposition \ref{produit1},
to $\widetilde{L}^{\rho}_{T}(B^{s_{1}}_{p,r})$ spaces. The indices $s_{1}$, $p$, $r$
behave just as in the stationary case whereas the time exponent $\rho$ behaves according to H\"older inequality.
\end{remarka}
In the sequel we will need of composition lemma in $\widetilde{L}^{\rho}_{T}(B^{s}_{p,r})$ spaces.
\begin{lemme}
\label{composition}
Let $s>0$, $(p,r)\in[1,+\infty]$ and $u\in \widetilde{L}^{\rho}_{T}(B^{s}_{p,r})\cap L^{\infty}_{T}(L^{\infty})$.
\begin{enumerate}
 \item Let $F\in W_{loc}^{[s]+2,\infty}(\T^{N})$ such that $F(0)=0$. Then $F(u)\in \widetilde{L}^{\rho}_{T}(B^{s}_{p,r})$. More precisely there exists a function $C$ depending only on $s$, $p$, $r$, $N$ and $F$ such that:
$$\|F(u)\|_{\widetilde{L}^{\rho}_{T}(B^{s}_{p,r})}\leq C(\|u\|_{L^{\infty}_{T}(L^{\infty})})\|u\|_{\widetilde{L}^{\rho}_{T}(B^{s}_{p,r})}.$$
\item Let $F\in W_{loc}^{[s]+3,\infty}(\T^{N})$ such that $F(0)=0$. Then $F(u)-F^{'}(0)u\in \widetilde{L}^{\rho}_{T}(B^{s}_{p,r})$. More precisely there exists a function $C$ depending only on $s$, $p$, $r$, $N$ and $F$ such that:
$$\|F(u)-F^{'}(0)u\|_{\widetilde{L}^{\rho}_{T}(B^{s}_{p,r})}\leq C(\|u\|_{L^{\infty}_{T}(L^{\infty})}\|u\|^{2}_{\widetilde{L}^{\rho}_{T}(B^{s}_{p,r})}.$$
\end{enumerate}
\end{lemme}
Here we recall a result of interpolation which explains the link
of the space $B^{s}_{p,1}$ with the space $B^{s}_{p,\infty}$, see
\cite{DFourier}.
\begin{proposition}
\label{interpolationlog}
There exists a constant $C$ such that for all $s\in\R$, $\e>0$ and
$1\leq p<+\infty$,
$$\|u\|_{\widetilde{L}_{T}^{\rho}(B^{s}_{p,1})}\leq C\frac{1+\e}{\e}\|u\|_{\widetilde{L}_{T}^{\rho}(B^{s}_{p,\infty})}
\biggl(1+\log\frac{\|u\|_{\widetilde{L}_{T}^{\rho}(B^{s+\e}_{p,\infty})}}
{\|u\|_{\widetilde{L}_{T}^{\rho}(B^{s}_{p,\infty})}}\biggl).$$ \label{5Yudov}
\end{proposition}
Now we give some result on the behavior of the Besov spaces via some pseudodifferential operator (see \cite{DFourier}).
\begin{definition}
Let $m\in\R$. A smooth function function $f:\T^{N}\rightarrow\R$ is said to be a ${\cal S}^{m}$ multiplier if for all muti-index $\alpha$, there exists a constant $C_{\alpha}$ such that:
$$\forall\xi\in\T^{N},\;\;|\p^{\alpha}f(\xi)|\leq C_{\alpha}(1+|\xi|)^{m-|\alpha|}.$$
\label{smoothf}
\end{definition}
\begin{proposition}
Let $m\in\R$ and $f$ be a ${\cal S}^{m}$ multiplier. Then for all $s\in\R$ and $1\leq p,r\leq+\infty$ the operator $f(D)$ is continuous from $B^{s}_{p,r}$ to $B^{s-m}_{p,r}$.
\label{singuliere}
\end{proposition}
We now focus on the mass equation associated to (\ref{1})
\begin{equation}
\begin{cases}
 \begin{aligned}
&\p_{t}q+u\cdot\n q+q h(q)=-q{\rm div}v_{1},\\
&q_{/ t=0}=q_{0}.
 \end{aligned}
\end{cases}
\label{21}
\end{equation}
where $h\in C^{\infty}$, $h(0)=0$ and $h^{'}\in W^{s,\infty}(\R,\R)$.  Here $v_{1}$ belongs in $\widetilde{L}^{1}(B^{\NN+\e}_{p,1})$ with $\e>0$ and $p\in[1,+\infty]$.
\begin{proposition}
Let 
$1\leq p\leq p_{1}\leq+\infty$, $p_{2}\in[1,+\infty]$ and $1\leq r\leq+\infty$, $1\leq r_{1}\leq+\infty$. Let assume that:
\begin{equation}
-N\min(\frac{1}{p_{1}},\frac{1}{p^{'}})<\sigma\leq\frac{N}{p_{2}},
 \label{3.21}
\end{equation}
with strict inequality if $r<+\infty$.
Assume that $q_{0}\in B^{\sigma}_{p,r}$, $\n v_{1}\in L^{1}(0,T;L^{\infty})$, ${\rm div}v_{1}\in \widetilde{L}^{1}_{T}(B^{\frac{N}{p_{1}}}_{p_{1},\infty})\cap L^{\infty}$  and that
$q\in\widetilde{L}^{\infty}_{T}(B^{\sigma}_{p,r})$
satisfies (\ref{21}).
There exists a constant $C$ depending only on $N$ such that for all
$t\in[0,T]$ and $m\in\mathbb{Z}$, we have:
\begin{equation}
\|q\|_{\widetilde{L}^{\infty}_{t}(B^{\sigma}_{p,r})}\leq e^{CV(t)}\|q_{0}\|_{B^{\sigma}_{p,r}}+e^{CV(t)}-1,
\label{22}
\end{equation}
with:
$$
\begin{cases}
\begin{aligned}
V(t)&=\int^{t}_{0}\big(\|\n u(\tau)\|_{B^{\frac{N}{p_{1}}}_{p_{1},\infty}\cap L^{\infty}}+\|{\rm div}v_{1}(\tau)\|_{B^{\frac{N}{p_{2}}}_{p_{2},\infty}\cap L^{\infty}}+\|q(\tau)\|^{\alpha+1}_{L^{\infty}}\big)d\tau\:\:\:\mbox{if}\;\;
\sigma<1+\frac{N}{p_{1}},\\
&=\int^{t}_{0}\big(\|\n u(\tau)\|_{B^{\frac{N}{p_{1}}}_{p_{1},\infty}\cap L^{\infty}}+\|{\rm div}v_{1}(\tau)\|_{B^{\frac{N}{p_{2}}}_{p_{2},\infty}\cap L^{\infty}}+\|q(\tau)\|^{\alpha+1}_{L^{\infty}}\big)d\tau\:\:\:\mbox{if}\;\;
\sigma\leq1+\frac{N}{p_{1}}\\
&\hspace{11cm}\;\;\mbox{and}\;\;r=1.
\end{aligned}
\end{cases}
$$
with $\alpha$ the smallest integer such that $\alpha\geq s$.
\label{transport}
 \end{proposition}
{\bf Proof:}
 Applying $\D_{l}$ to (\ref{21}) yields:
$$\p_{t}\D_{l}q+u\cdot\n\D_{l}q+\D_{l}q=R_{l}-\D_{l}(q{\rm div}v_{1})-\D_{l}(q h(q))\;\;\;\mbox{with}\;\;R_{l}=[u\cdot\n,\D_{l}]q.$$
Multiplyng by $\D_{l}q|\D_{l}q|^{p-2}$ and performing a time integration, we easily get:
$$
\begin{aligned}
&\|\D_{l}q(t)\|_{L^{p}}d\tau\lesssim\|\D_{l}q_{0}\|_{L^{p}}+\int^{t}_{0}\big(\|R_{l}\|_{L^{p}}
+\|{\rm div}u\|_{L^{\infty}}\|\D_{l}q\|_{L^{p}}\\
&\hspace{7cm}+\|\D_{l}(q{\rm div}v_{1})\|_{L^{p}}+\|\D_{l}(q h(q))\|_{L^{p}}\big)d\tau.
\end{aligned}
$$
By paraproduct, there exists a constant $C$ and a positive sequence $(c_{l})\in l^{r}$ such that:
$$\|\D_{l}(q{\rm div}v_{1})\|_{L^{p}}\leq C c_{l}2^{-l\sigma}\|q\|_{B^{\sigma}_{p,r}}\|{\rm div}v_{1}\|_{B^{\frac{N}{p_{2}}}_{p_{2},\infty}\cap L^{\infty}}.$$
Similarly by lemma \ref{composition} we have:
$$\|\D_{l}(qh(q))\|_{L^{p}}\leq C c_{l}2^{-q\sigma}\|q\|_{B^{\sigma}_{p,r}}\|q\|^{\alpha+1}_{L^{\infty}}.$$
Next the term $\|R_{l}\|_{L^{p}}$ may be bounded according to the inequality (2.53) p 110 of \cite{BCD}:
$$
\begin{aligned}
\|\big(2^{l\sigma}\|R_{l}\|_{L^{p}}\big)_{l}\|_{l^{r}}&\leq C\|\n u \|_{B^{\frac{N}{p_{1}}_{p_{1},r}}}\|q\|_{B^{\sigma}_{p,r}}.
\end{aligned}
$$
We end up with multiplying the previous inequality by $2^{l(\NN+\e)}$ and summing up on $\mathbb{Z}$:
$$\|q(t)\|_{B^{\sigma}_{p,r}}\leq \|q_{0}\|_{B^{\sigma}_{p,r}}
+\int^{t}_{0}C V^{'}\|q\|_{B^{\sigma}_{p,r}}d\tau+\int^{t}_{0}CV^{'}d\tau.$$
Gronwall lemma yields inequality (\ref{22}).
\hfill{$\Box$}
\section{A priori bounds on the density}
\label{section3}
In this section we make a formal analysis on the partial differential equations of (\ref{1}) and begin by classical energy estimates. Multiplying the equation of conservation of momentum by $u$, we obtain
 \begin{equation}
\begin{aligned}
&\int_{\T^{N}}(\frac{1}{2}\rho
|u|^{2}(t,x)+\Pi(\rho)(t,x)dx+\int_{0}^{t}\int_{\T^{N}}(\mu
D(u):D(u)(s,x)\\
&\hspace{2cm}+(\lambda+\mu)|{\rm div} u|^{2}(s,x))dsdx
\leq\int_{\T^{N}}\big(\frac{|m_{0}|^{2}}{2\rho}(x)+\Pi(\rho_{0})(x)]\big)dx,
\end{aligned}
\label{17}
\end{equation}
where $\Pi$ is defined by
\begin{equation}
 \Pi(s)=s\biggl(\int^{s}_{0}\frac{P(z)}{z^{2}}dz\biggl),
\end{equation}
It follows classicaly that we have the following bounds
\begin{equation}
 \begin{cases}
  \begin{aligned}
&\rho\in L^{\infty}(0,\infty;L^{\gamma}_{2}),\\
&\sqrt{\rho}u\in L^{\infty}(0,\infty;L^{2}),\\
&\n u\in L^{2}(0,\infty;L^{2})^{N^{2}}.
  \end{aligned}
\end{cases}
\label{19}
\end{equation}
\subsection{Bound on $\log\rho$}
Let us emphasize that one of the main ingredient of the proof of the first part of theorem \ref{2} is a partial differential equation derived from (\ref{1}) involving $\log\rho$. It was introduced by P-L Lions in \cite{13} to prove global existence of weak solutions of (\ref{1}) in a particular case and it is one of the key of the proof in the paper of B. Desjardins in \cite{CCK5}. Letting formally $(\D)^{-1}{\rm div}$ operate on the equation of conservation of momentum, we obtain
\begin{equation}
(2\mu+\lambda){\rm div}u-P(\rho)+\int_{\T^{N}}P(\rho)dx=\p_{t}\D^{-1}{\rm
div}(\rho\,u)+R_{i}R_{j}(\rho\,u_{i}u_{j})
\label{20}
\end{equation}
where $\D^{-1}$ denotes the inverse Laplacian with zero mean value on ${\cal T}^{N}$
and $R_{i}$ the usual Riesz transform. Let us observe that the equation of mass yields
\begin{equation}
\p_{t}\log\rho+u.\n\log\rho+{\rm div}u=0
 \label{g21}
\end{equation}
Let us define $F$ and $G$ by the following expression:
$$
\begin{aligned}
&F=(2\mu+\lambda)(\log\rho+\D^{-1}{\rm
div}(\rho\,u)),\\
&G=(2\mu+\lambda){\rm
div}u-P(\rho).
\end{aligned}
$$
Moreover we shall denote respectively by $P$ and $Q$ the projection on the space of divergence-free and curl-free vector fields.
Combining (\ref{20}) and (\ref{21}), we obtain
\begin{equation}
\p_{t}F+u\cdot \n F+P(\rho)-\int_{\T^{N}}P(\rho)dx=[u_{j},R_{i}R_{j}](\rho u_{i}).
 \label{24}
\end{equation}
Next we define the Lagrangian flow $X$ of $u$ by
\begin{equation}
\begin{cases}
 \begin{aligned}
 &\p_{t}X(t,s,x)=u(t,X(t,s,x)),\\
&X_{/ t=s}=x,
 \end{aligned}
\end{cases}
\label{25}
\end{equation}
and derive the following identity
\begin{equation}
F(t,X(t,0,x))=F_{0}(x)-\int^{t}_{0}P(\rho(s,\cdot))dxds+\int^{t}_{0}([u_{j},R_{i}R_{j}](\rho u_{i})(s,X(s,0,x))ds.
 \label{26}
\end{equation}
Using the fact that $\rho(\cdot)\geq 0$, we obtain
\begin{equation}
F(t,x)\leq F_{0}(X(0,t,x))+\int^{t}_{0}P(\rho(s,\cdot))dxds+\int^{t}_{0}([u_{j},R_{i}R_{j}](\rho u_{i})(s,X(s,t,x))ds.
 \label{27}
\end{equation}
It follows that
\begin{equation}
\begin{aligned}
\log(\rho(t,x))\leq&\log(\|\rho_{0}\|_{L^{\infty}})+C\|(\D)^{-1}{\rm div} m_{0}\|_{L^{\infty}}+C\|(\D)^{-1}{\rm div}(\rho u)(t,\cdot)\|_{L^{\infty}}\\
&+C_{0}t+C \int^{t}_{0}\|[u_{j},R_{i}R_{j}](\rho u_{i})(s,\cdot)\|_{L^{\infty}}ds,
\end{aligned}
 \label{g28}
\end{equation}
where $C_{0}$ depends of the initial data. In view of the usual Sobolev embedding inequalities, we obtain
\begin{equation}
\begin{aligned}
\log(\rho(t,x))\leq&\log(\|\rho_{0}\|_{L^{\infty}})+C\|(\D)^{-1}{\rm div} m_{0}\|_{L^{\infty}}+C\|(\D)^{-1}{\rm div}(\rho u)\|_{L^{\infty}}\\
&+C_{0}t+C \int^{t}_{0}\|[u_{j},R_{i}R_{j}](\rho u_{i})(s,\cdot)\|_{B^{1}_{N+\e,1}}ds,
\end{aligned}
\label{29}
\end{equation}
with $\e>0$. Let us now remark that we have
\begin{equation}
\|\n u\|_{L^{N+\e}}\leq C(\|{\rm curl} u\|_{L^{N+\e}}+\|G\|_{L^{N+\e}}+\|P(\rho)\|_{L^{N+\e}}).
 \label{31}
\end{equation}
In view of R. Coifman, P-.L. Lions and S. Semmes \cite{1}, the following map
\begin{equation}
\begin{aligned}
&W^{1,r_{1}}(\T^{N})^{N}\times L^{r_{2}}(\T^{N})\rightarrow W^{1,r_{3}}(\T^{N})^{N}\\
&\hspace{3cm}(a,b)\rightarrow[a_{j},R_{i}R_{j}]b_{i}
\end{aligned}
\label{a33}
\end{equation}
is continuous for any $N\geq2$ as soon as $\frac{1}{r_{3}}=\frac{1}{r_{1}}+\frac{1}{r_{2}}$. Hence we have the following estimates for $N=3$:
\begin{equation}
\|[u_{j},R_{i}R_{j}](\rho u_{i})(s,\cdot)\|_{L^{1}(B^{1}_{N,1})}\leq C\|\n u\|_{L^{2}(L^{6})}\|\rho u\|_{L^{2}(L^{6+\e})},
 \label{33}
\end{equation}
with $\e>0$.
We obtain finally:
\begin{equation}
\begin{aligned}
\log(\rho(t,x))\leq&\log(\|\rho_{0}\|_{L^{\infty}})+C\|(\D)^{-1}{\rm div} m_{0}\|_{L^{\infty}}+C\|(\D)^{-1}{\rm div}(\rho u)\|_{L^{\infty}}\\
&\hspace{4,5cm}+C_{0}t+C \|\n u\|_{L^{2}_{t}(L^{6})}\|\rho u\|_{L^{2}_{t}(L^{6+\e})},
\end{aligned}
\label{29imp}
\end{equation}
\section{A priori estimates for the velocity}
\label{section4}
\subsection{Gain of integrability of the velocity $u$}
We want here derive estimate of integrability on the velocity $u$. This idea has been successively used in different papers,
we refer in particularly to \cite{5H2} and \cite{5MV,5MV2}. To do it, we multiply the momentum equation by $u|u|^{p_{1}-2}$ and we get after integration by part:
$$
\begin{aligned}
&\frac{1}{p_{1}}\int_{\T^{N}}\rho|u|^{p_{1}}(t,x)dx+\mu\int^{t}_{0}\int_{\T^{N}}\big(|u|^{p_{1}-2}|\n u|^{2}(t,x)+\frac{p_{1}-2}{4}
|u|^{p_{1}-4}|\n|u|^{2}|^{2}(t,x)\big)dxdt\\
&\hspace{1,8cm}+\lambda\int^{t}_{0}\int_{\T^{N}}\big(({\rm div}u)^{2}|u|^{p_{1}-2}(t,x)+\frac{p_{1}-2}{2}
{\rm div}u\sum_{i}u_{i}\p_{i}|u|^{2}|u|^{p_{1}-4}(t,x)\big)dtdx\\
&\hspace{2,5cm}-\int^{t}_{0}\int_{\T^{N}} P(\rho)\big({\rm div}u|u|^{p_{1}-2}+(p_{1}-2)
\sum_{i,k}u_{i}u_{k}\p_{i}u_{k}|u|^{p_{1}-4}\big)(t,x)dtdx\\
&\hspace{11cm}\leq\int_{\T^{N}}\rho_{0}|u_{0}|^{p_{1}}dx.
\end{aligned}
$$
We have then by Young's inequality:
$$
\begin{aligned}
&\frac{\lambda(p_{1}-2)}{2}\int^{t}_{0}\int_{\T^{N}}
{\rm div}u\sum_{i}u_{i}\p_{i}|u|^{2}|u|^{p_{1}-4}(t,x)\big)dtdx\\
&=\lambda\frac{p_{1}-2}{2}\int^{t}_{0}\int_{\T^{N}}
{\rm div}u\,u\cdot\n (|u|^{2})|u|^{p_{1}-4}(t,x)dtdx\leq\\
&\lambda \frac{p_{1}-2}{2}(\frac{\eta}{2}\int^{t}_{0}\int_{\T^{N}}
|{\rm div}u|^{2}|u|^{p_{1}-2}dtdx+\frac{2}{\eta}\int^{t}_{0}\int_{\T^{N}}|\n |u|^{2}|^{2}|u|^{p_{1}-4}(t,x)dtdx)\\
\end{aligned}
$$
If we choose:
$$\lambda \frac{\eta(p_{1}-2)\lambda}{4}=s\mu+\lambda,$$
for some $s\in(0,\frac{1}{N})$,by the fact that $({\rm div}u)^{2}\leq N |\n u|^{2}$ we therefore obtain:
$$
\begin{aligned}
&\frac{1}{p_{1}}\int_{\T^{N}}\rho|u|^{p_{1}}(t,x)dx+A_{s}\int^{t}_{0}\int_{\T^{N}}|u|^{p_{1}-2}|\n u|^{2}(t,x)dtdx\\
&\hspace{1cm}+B_{s}\int^{t}_{0}\int_{\T^{N}}
|u|^{p_{1}-4}|\n|u|^{2}|^{2}(t,x)dxdt\leq
\int^{t}_{0}\int_{\T^{N}} P(\rho)\big({\rm div}u|u|^{p_{1}-2}\\
&\hspace{3cm}+\frac{p_{1}-2}{2}
u\cdot\n(|u|^{2})|u|^{p_{1}-4}\big)(t,x)dtdx+\int_{\T^{N}}\rho_{0}|u_{0}|^{p_{1}}dx,
\end{aligned}
$$
with $A_{s}=\mu(1-sN)$ and $B(s)=\frac{p_{1}-2}{4}\mu-\frac{(p_{1}-2)^{2}\lambda^{2}}{16(s\mu+\lambda)}$.
By Young inequality, we get again:
$$
\begin{aligned}
&\frac{1}{p_{1}}\int_{\T^{N}}\rho|u|^{p_{1}}(t,x)dx+A_{s}\int^{t}_{0}\int_{\T^{N}}|u|^{p_{1}-2}|\n u|^{2}(t,x)dtdx\\
&+B_{s}\int^{t}_{0}\int_{\T^{N}}
|u|^{p_{1}-4}|\n|u|^{2}|^{2}(t,x)dxdt\leq
C_{\e}\int^{t}_{0}\int_{\T^{N}} P(\rho)^{2}|u|^{p_{1}-2}dtdx+\int_{\T^{N}}\rho_{0}|u_{0}|^{p_{1}}dx,
\end{aligned}
$$
with $C_{\e}$ enough big. Now we want use the fact that $\n |u|^{\frac{p_{1}}{2}}\in L^{2}_{t}(L^{2})$, which implies that when $N=3$ $u\in L^{p_{1}}_{t}(L^{3p_{1}})$. More precisely we have:
$$\|u\|^{\frac{p_{1}}{2}}_{L^{p_{1}}_{t}(L^{3p_{1}})}\leq C(\|\n(|u|^{\frac{p_{1}}{2}})\|_{L^{2}(L^{2})}+\bar{u}_{\frac{p_{1}}{2}}),$$
where $\bar{u}_{\frac{p_{1}}{2}})$  is the average of $|u|^{\frac{p_{1}}{2}}$. We have then by H\"older's inequalities with $\frac{p_{1}-2}{3p_{1}}+\frac{2(p_{1}+1)}{3p_{1}}=1$ and $\frac{p_{1}-2}{p_{1}}+\frac{2}{p_{1}}=1$:
$$
\begin{aligned}
 |\int^{t}_{0}\int_{\T^{N}} P(\rho)^{2}|u|^{p_{1}-2}dtdx|&\leq
\|P(\rho)^{2}\|_{L_{t}^{\frac{p_{1}}{2}}(L^{\frac{3p_{1}}{2(p_{1}+1)}})}
\||u|^{p_{1}-2}\|_{L_{t}^{\frac{p_{1}}{p_{1}-2}}(L^{\frac{3p_{1}}{p_{1}-2}})},\\
&\leq\|P(\rho)\|_{L_{t}^{p_{1}}(L^{\frac{3p_{1}}{p_{1}+1}})}^{2}\|u\|_{L^{p_{1}}_{t}(L^{3p_{1}})}^{p_{1}-2},\\
&\leq C\|P(\rho)\|_{L_{t}^{p_{1}}(L^{\frac{3p_{1}}{p_{1}+1}})}^{2}(
\|\n( |u|^{\frac{p_{1}}{2}})+\widetilde{|u|^{\frac{p_{1}}{2}}}(\cdot)\|_{L_{t}^{2}(L^{2})})^{2-\frac{4}{p_{1}}}.
\end{aligned}
$$
Remarking that $\int_{\T^{N}}\rho_{0}dx=M\ne 0$, we can write as $\gamma\geq\frac{6}{5}$:
$$
\begin{aligned}
&\widetilde{|u|^{\frac{p_{1}}{2}}}(s)\leq\frac{1}{M}(\|\rho(s,\cdot)\|_{L^{\gamma}}\|\n |u|^{\frac{p_{1}}{2}}(s,\cdot)\|_{L^{2}}+\int_{\T^{N}}\rho(s,x)|u(s,x)|^{\frac{p_{1}}{2}}dx),\\
&\big(\int^{t}_{0}\widetilde{|u|^{\frac{p_{1}}{2}}}^{2}(s)ds\big)^{\frac{1}{2}}\leq
\frac{1}{M}(\|\rho\|_{L^{\infty}_{t}(L^{\gamma})}\|\n |u|^{\frac{p_{1}}{2}}\|_{L^{2}_{t}(L^{2})}+(Mt)^{\frac{1}{2}}\|\rho^{\frac{1}{p_{1}}}u\|_{L^{\infty}_{t}(L^{p_{1}})}^{\frac{p_{1}}{2}}),\\
\end{aligned}
$$
We obtain then:
$$
\begin{aligned}
&\|P(\rho)\|_{L_{t}^{p_{1}}(L^{\frac{3p_{1}}{p_{1}+1}})}^{2}(
\|\n( |u|^{\frac{p_{1}}{2}})+\widetilde{|u|^{\frac{p_{1}}{2}}}(\cdot)\|_{L_{t}^{2}(L^{2})})^{2-\frac{4}{p_{1}}}\leq\\
&C_{1}\|P(\rho)\|_{L_{t}^{p_{1}}(L^{\frac{3p_{1}}{p_{1}+1}})}^{2}\big(
\|\n( |u|^{\frac{p_{1}}{2}})\|_{L_{t}^{2}(L^{2})}^{\frac{2p_{1}-4}{p_{1}}}+
(Mt)^{\frac{p_{1}-2}{p_{1}}}\|\rho^{\frac{1}{p_{1}}}u\|_{L^{\infty}_{t}(L^{p_{1}})}^{p_{1}-2})
\end{aligned}
$$
By a standard application of Young inequality ($\frac{2p_{1}-4}{2p_{1}}+\frac{4}{2p_{1}}=1$), we obtain that:
\begin{equation}
\begin{aligned}
&\frac{1}{p_{1}}\int_{\T^{N}}\rho|u|^{p_{1}}(t,x)dx+A_{s}\int^{t}_{0}|u|^{p_{1}-2}|\n u|^{2}(t,x)dtdx\\
&+B_{s}\int^{t}_{0}
|u|^{p_{1}-4}|\n|u|^{2}|^{2}(t,x)dxdt\leq C_{\e,t}^{2}\|P(\rho)\|_{L_{t}^{p_{1}}(L^{\frac{3p_{1}}{p_{1}+1}})}^{2p_{1}}+\frac{1}{p_{1}}\int_{\T^{N}}\rho_{0}|u_{0}|^{p_{1}}dx,\\
\end{aligned}
\label{againvitesse}
\end{equation}
where $C_{\e,t}$ is big enough and depend of the time $t$.
\subsection{Gain of derivatives on the velocity $u$}
\label{gainvitesse}
In this section we deal with the case $N=3$. The case $N=2$ follows the same lines. In the sequel we will follow the procedure
developped in \cite{CCK5} and \cite{H2} to get some energy inequalities. The main idea compared with the results in \cite{CCK2, CCK3,5CK2} is to obtain energy inequalities which depends only of the control on $\rho\in L^{\infty}$. It implies that we have to be carreful to not introduce some derivatives on the density which means to ``kill'' the coupling between velocity and pressure. Multiplying first the equation of conservation of momentum by $f(t)\p_{t}u$ with $f(t)=\min(1,t)$ 
and integrating over $(0,T)\times\T^{N}$, we deduce that:
\begin{equation}
 \begin{aligned}
  &\int^{t}_{0}\int_{\T^{N}}f(t)\rho|\p_{t}u|^{2}dxds+\frac{1}{2}\int_{\T^{N}}f(t)\big(\mu|\n u(t,x)|^{2}+(\lambda+\mu)|{\rm div}u|^{2}(t,x)\big)dx\\
&+\int^{t}_{0}\int_{\T^{N}}\n P(\rho)\cdot f(t)\p_{t}u dxds\leq \int^{t}_{0}\int_{\T^{N}}\frac{f^{'}(t)}{2}\big(\mu|\n u(t,x)|^{2}+(\lambda+\mu)|{\rm div}u|^{2}(t,x)\big)dxds\\
&\hspace{2cm}+\int^{t}_{0}\|\sqrt{f(t)\rho}\p_{t}u\|_{L^{2}(\T^{N})}
(\|\sqrt{\rho}(u\cdot\n)u\|_{L^{2}(\T^{N})}+\|\sqrt{\rho}f\|_{L^{2}(\T^{N})})ds.
 \end{aligned}
\label{a74}
\end{equation}
Next we use the equation of mass conservation to write:
$$
\begin{aligned}
& \int^{t}_{0}\int_{\T^{N}}\n P(\rho)\cdot f(t)\p_{t}u dxds=-\p_{t}\int_{\T^{N}}f(t)P(\rho){\rm div}u dx+\int_{\T^{N}}\p_{t}(f(t)P(\rho)){\rm div}u dx\\
&=-\p_{t}\int_{\T^{N}}f(t)P(\rho){\rm div}u dx-\int_{\T^{N}}f(t)\big[{\rm div}(P(\rho)u){\rm div}u+(\rho P^{'}(\rho)-P(\rho)){\rm div}u\big] dx\\
&\hspace{10cm}-\int_{\T^{N}}f^{'}(t)P(\rho){\rm div}u dx,\\
&=-\p_{t}\int_{\T^{N}}f(t)P(\rho){\rm div}u dx+\frac{1}{2\mu+\lambda}\int_{\T^{N}}f(t)P(\rho)u\cdot\n (G+P(\rho)) dx\\
&-\frac{1}{(2\mu+\lambda)^{2}}\int_{\T^{N}}f(t)(\rho P^{'}(\rho)-P(\rho))\big(G^{2}-P(\rho)^{2}+2(\lambda+2\mu)P(\rho){\rm div}u\big)dx \\
&\hspace{9cm}-\int_{\T^{N}}f^{'}(t)P(\rho){\rm div}u dx,
\end{aligned}
$$
If we define $\Pi_{f}(s)=s(\int^{s}_{0}\frac{f(z)}{z^{2}}dz)$, we have then by mass equation:
$$\p_{t}\Pi_{f}(\rho)+{\rm dib}(\Pi_{f}(\rho)u)+f(\rho){\rm div}u=0.$$
By this fact we obtain that:
$$P(\rho)\big(u\cdot\n P(\rho)-2(\rho P^{'}(\rho)-P(\rho)){\rm div}u\big)=P(\rho)(-\p_{t}(P(\rho))-3P^{'}(\rho)\rho{\rm div}u+2P(\rho){\rm div}u).$$
and:
$$
\begin{aligned}
&\int_{\T^{N}}\n P(\rho)\cdot f(t)\p_{t}u dx=-\p_{t}\int_{\T^{N}}f(t)P(\rho){\rm div}u dx+\frac{1}{\lambda+2\mu}\p_{t}\int_{\T^{N}}f(t)k(\rho)dx\\
&+\frac{1}{2\mu+\lambda}\int_{\T^{N}}f(t)P(\rho)u\cdot\n G dx+\frac{1}{(2\mu+\lambda)^{2}}\int_{\T^{N}}f(t)P^{2}(\rho)(\rho P^{'}(\rho)-P(\rho)) dx\\
&\hspace{1,8cm}-\frac{1}{(2\mu+\lambda)^{2}}\int_{\T^{N}}f(t)G^{2}(\rho P^{'}(\rho)-P(\rho)) dx-\frac{1}{\lambda+2\mu}\int_{\T^{N}}f^{'}(t)k(\rho)dx\\
&\hspace{9cm}-\int_{\T^{N}}f^{'}(t)P(\rho){\rm div}u dx.
\end{aligned}
$$
where $k(s)=P(s)^{2}-P(\bar{\rho})^{2}+s(\int^{s}_{\bar{\rho}}\frac{P(s)}{s^{2}}ds-\frac{P(\bar{\rho})}{\bar{\rho}})$.
Inserting the above inequality in (\ref{a74}) and by Young's inequality, we obtain:
\begin{equation}
 \begin{aligned}
  \int^{t}_{0}\int_{\T^{N}}f(t)\rho|\p_{t}u|^{2}dxds+\frac{1}{2}\int_{\T^{N}}f(t)\big(\mu|\n u(t,x)|^{2}+(\lambda+\mu)({\rm div}u(t,x))^{2}\big)dx&\\
+\frac{1}{(2\mu+\lambda)^{2}}\int^{t}_{0}\int_{\T^{N}}f(t)P(\rho)^{2}(\rho P^{'}(\rho)-P(\rho))dx ds&\\
+\frac{1}{\lambda+2\mu}\int_{\T^{N}}f(t)k(\rho(t,x))dx\leq C+\int_{\T^{N}}f(t)P(\rho(t,x)){\rm div}u(t,x)dx&\\
+\frac{1}{\lambda+2\mu}\int^{t}_{0}\int_{\T^{N}}f^{'}(t)k(\rho)dx+\int^{t}_{0}\int_{\T^{N}}f^{'}(t)P(\rho){\rm div}u dx&\\
+C\int^{t}_{0}\int_{\T^{N}}f(t)(|\rho P^{'}(\rho)-P(\rho)|G^{2}+|P(\rho)u||\n G|+|\sqrt{\rho}u\cdot\n u|^{2}&\\
+|\sqrt{\rho}f|^{2})dxds.&
 \end{aligned}
\label{79}
\end{equation}
In the sequel we set:
$$
\begin{aligned}
&A(t)=\int^{t}_{0}\int_{\T^{N}}f(t)\rho|\p_{t}u|^{2}dxds+\frac{1}{2}\int_{\T^{N}}f(t)\big(\mu|\n u(t,x)|^{2}+(\lambda+\mu)({\rm div}u(t,x))^{2}\big)dx\\
&+\frac{1}{(2\mu+\lambda)^{2}}\int^{t}_{0}\int_{\T^{N}}f(s)P(\rho)^{2}(\rho P^{'}(\rho)-P(\rho))dx ds+\frac{1}{\lambda+2\mu}\int_{\T^{N}}f(t)k(\rho(t,x))dx
\end{aligned}
$$
We obtain finally:
\begin{equation}
 \begin{aligned}
&A(t)\leq C+C_{t}\|\rho\|_{L^{\infty}}+C\int^{t}_{0}\int_{\T^{N}}f(s)(|\rho P^{'}(\rho)-P(\rho)|G^{2}+|P(\rho)u||\n G|\\
&\hspace{8cm}+|\sqrt{\rho}u\cdot\n u|^{2}+|\sqrt{\rho}f|^{2})dxds.\\
&\leq C+C\int^{t}_{0}(\|\rho P^{'}(\rho)-P(\rho)\|_{L^{\infty}}\|\sqrt{f(s)}G\|_{L^{2}}^{2}+f(s)\|g(\rho)(s,\cdot)\|_{L^{\infty}}
\|\sqrt{\rho}u\|_{L^{2}}\\
&\hspace{2cm}\times\|\n G\|_{L^{2}}+\|\sqrt{\rho}u\|_{L^{4}}^{2}\|\sqrt{f(s)}\n u\|_{L^{4}}^{2}+\|\rho\|_{L^{\infty}}\|\sqrt{f(s)}f\|_{L^{2}}^{2})dxds,\\
&\leq C+C\int^{t}_{0}(\|h(\rho(s,\cdot)\|_{L^{\infty}}\|\sqrt{f(s)}\n u\|_{L^{2}}^{2}+f(s)\|i(\rho(s,\cdot)\|_{L^{\infty}}+f(s)\|\n G\|_{L^{2}}
\\
&\hspace{1cm}\times \|g(\rho(s,\cdot))\|_{L^{\infty}}+\|\sqrt{\rho}u\|_{L^{4}}^{2}\|\sqrt{f(s)}\n u\|_{L^{4}}^{2}+\|\rho\|_{L^{\infty}}\|\sqrt{f(s)}f\|_{L^{2}}^{2})dxds,
 \end{aligned}
\label{a81}
\end{equation}
where $g(\rho)=\frac{P(\rho)}{\sqrt{\rho}}$, $h(s)=|s P^{'}(s)-P(s)|$ and $i(s)=h(s)P(s)^{2}$.
\subsubsection*{Estimates on $Pu$ and $G$}
We want now to obtains bounds on ${\cal P}u$ and $g$, assuming that $\rho$ is a priori bounded in $L^{\infty}(\T^{N})$.
Indeed we wnat show that the control of $A(t)$ in (\ref{a81}) depend only of a control on $\|\rho\|_{L^{\infty}}$.\\
Next we use once more the equation of conservation of momentum to write:
\begin{equation}
\mu\D u={\cal P}(\rho\p_{t}u)+{\cal P}(\rho u\cdot\n u)-{\cal P}(\rho f),
 \label{83}
\end{equation}
\begin{equation}
\n G= {\cal Q}(\rho\p_{t}u)+{\cal Q}(\rho u\cdot\n u)-{\cal Q}(\rho f).
 \label{84}
\end{equation}
Therefore we have:
\begin{equation}
\begin{aligned}
&\|\n G\|_{L^{2}}+\|\D{\cal P}u\|_{L^{2}}\leq C\|\rho(s,\cdot)\|^{\frac{1}{2}}_{L^{\infty}}\big(\|\sqrt{\rho}\p_{s}u(s\cdot)\|_{L^{2}}
+\|\sqrt{\rho}u\cdot\n u(s,\cdot)\|_{L^{2}}\\
&\hspace{5cm}+\|\rho(s,\cdot)\|^{\frac{1}{2}}_{L^{\infty}}\|f(s,\cdot)\|_{L^{2}}\big).
\end{aligned}
 \label{86}
\end{equation}
\subsubsection*{The case $N=3$}
For simplicity, we will treat only the case of the dimension $3$. We recall that for all $1<p<+\infty$:
$$\|\n u\|_{L^{p}}\leq \|\n {\cal P}u\|_{L^{p}}+\|R G\|_{L^{p}}+\|R(P(\rho))\|_{L^{p}}.$$
We want now to recall the Gagliardo-Nirenberg's theorem:
$$\forall f\in H^{1}(\T^{N})\;\;\mbox{such that}\;\;\int_{\T^{N}}f dx=0,\;\;\|f\|^{2}_{L^{4}(\T^{N})}\leq C\|f\|^{\frac{1}{2}}_{L^{2}(\T^{N})}
\|\n f\|_{L^{2}(\T^{N})}^{\frac{3}{2}}.$$
We deduce that from Gagliardo-Nirenberg's inequality and Young's inequalities:
\begin{equation}
\begin{aligned}
&\|\sqrt{\rho}u\|^{2}_{L^{4}}\|\sqrt{f(s)}\n u\|^{2}_{L^{4}}\leq Cf(s)\|\sqrt{\rho}u\|^{2}_{L^{4}}(\|R(P(\rho))\|^{2}_{L^{4}}
+\|\n {\cal P}u\|_{L^{4}}+\|R G\|_{L^{4}}\big)\\
&\leq C\|\sqrt{\rho}u\|^{2}_{L^{4}}\big(f(s)\|P(\rho)\|^{2}_{L^{4}}+\big(\sqrt{f(s)}(\|\n u\|_{L^{4}}+\|P(\rho)\|_{L^{2}}))^{\frac{1}{2}}\\
&\hspace{6cm}\times(\sqrt{f(s)}(\|\D {\cal P}u\|_{L^{2}}+\|\n G\|_{L^{2}})\big)^{\frac{3}{2}}\big)\\
&\leq C\big(f(s)\|\sqrt{\rho}u\|^{2}_{L^{4}}\|P(\rho)\|^{2}_{L^{4}}+\|\rho(s,\cdot)\|^{\frac{3}{4}}_{L^{\infty}}\|\sqrt{\rho}u\|^{2}_{L^{4}}\big(\sqrt{f(s)}(\|\n u\|_{L^{2}}+\|P(\rho)\|_{L^{2}}))^{\frac{1}{2}}\\
&\times(\sqrt{f(s)}(\|\sqrt{\rho}\p_{s}u(s\cdot)\|_{L^{2}}
+\|\sqrt{\rho}u\cdot\n u(s,\cdot)\|_{L^{2}}+\|\rho(s,\cdot)\|^{\frac{1}{2}}_{L^{\infty}}\|f(s,\cdot)\|_{L^{2}}\big)^{\frac{3}{2}}\big)\\
&\leq C\big(f(s)\|\sqrt{\rho}u\|^{2}_{L^{4}}\|P(\rho)\|^{2}_{L^{4}}+\frac{1}{\e}f(s)(\|\n u\|^{2}_{L^{2}}+\|P(\rho)\|^{2}_{L^{2}})\|\rho(s,\cdot)\|^{3}_{L^{\infty}}\|\sqrt{\rho}u\|^{8}_{L^{4}}\\
&+\e f(s)(\|\sqrt{\rho}\p_{s}u(s\cdot)\|^{2}_{L^{2}}+\|\sqrt{\rho}u\cdot\n u(s,\cdot)\|^{2}_{L^{2}}+\|\rho(s,\cdot)\|_{L^{\infty}}\|f(s,\cdot)\|^{2}_{L^{2}})\big).\\
\end{aligned}
 \label{99}
\end{equation}
Hence we obtain by Young inequality:
\begin{equation}
\begin{aligned}
&f(s)\|g(\rho(s,\cdot)\|_{L^{\infty}}\|\n G(s,\cdot)\|_{L^{2}}\leq \frac{C}{\e}\|g(\rho(s,\cdot)\|^{2}_{L^{\infty}}\|\rho(s,\cdot)\|_{L^{\infty}}f(s)\\
&+\e(\|\sqrt{f(s)\rho}\p_{s}u(s\cdot)\|^{2}_{L^{2}}
+\|\sqrt{\rho f(s)}u\cdot\n u(s,\cdot)\|^{2}_{L^{2}}+f(s)\|\rho(s,\cdot)\|_{L^{\infty}}\|f(s,\cdot)\|^{2}_{L^{2}}).
\end{aligned}
\label{100}
\end{equation}
By adding (\ref{99}) and (\ref{100}), we obtain:
\begin{equation}
\begin{aligned}
&\|\sqrt{\rho}u\|^{2}_{L^{4}} \|\sqrt{f(s)}\n u\|^{2}_{L^{4}} +f(s)\|g(\rho(s,\cdot)\|_{L^{\infty}}\|\n G(s,\cdot)\|_{L^{2}}\leq\\
&C\big(f(s)\|\sqrt{\rho}u\|^{2}_{L^{4}} \|P(\rho)\|^{2}_{L^{4}}+\frac{1}{\e}f(s)(\|\n u\|^{2}_{L^{2}}+\|P(\rho)\|^{2}_{L^{2}})\|\rho(s,\cdot)\|^{3}_{L^{\infty}}\|\sqrt{\rho}u\|^{8}_{L^{4}}\\
&+\e f(s)(\|\sqrt{\rho}\p_{s}u(s\cdot)\|^{2}_{L^{2}}+\|\rho(s,\cdot)\|_{L^{\infty}}\|f(s,\cdot)\|^{2}_{L^{2}})\big)+\frac{C}{\e}f(s)\|g(\rho(s,\cdot)\|^{2}_{L^{\infty}}\\
&\hspace{10cm}\times\|\rho(s,\cdot)\|_{L^{\infty}}.
\end{aligned}
\label{aa81}
\end{equation}
Therefore we have by using inequality (\ref{againvitesse}), (\ref{a81}) and (\ref{aa81}):
\begin{equation}
 \begin{aligned}
&A(t)\leq C+C\int^{t}_{0}\phi(\|\rho(s,\cdot)\|_{L^{\infty}})f(s)\big((1+\|\n u\|_{L^{2}}^{2}+\|P(\rho(s,\cdot))\|^{2}_{L^{2}})\\
&\hspace{10cm}+\|f(s,\cdot)\|_{L^{2}}^{2}\big)ds,
 \end{aligned}
\label{81}
\end{equation}
where $C$ depends of the time $t$ here.
Gronwall's lemma provides the following bound:
\begin{equation}
A(t)\leq \exp(C\exp(\int^{t}_{0}\phi(\|\rho(s,\cdot)\|_{L^{\infty}})ds)),
\label{5imp}
\end{equation}
where $\phi\in C^{0}(\R_{+},\R_{+}^{*})\cap C^{1}(0,\infty)$ such that $\phi(s)\geq \e_{0}s$ for some positive $s$.
\subsubsection*{Control of $\sup_{0<t\leq T}f^{2}(t)\int|\dot{u}|^{2}(t,x)dx+\int\int f^{2}(s)|\n\dot{u}|^{2}dxds$}
In the sequel, we want obtain estimate on $\n u$ in $L^{1}_{T}(BMO)$, that's why we need of additional regularity estimates.
We derive then estimates for the terms $f(t)^{2}\int_{\T^{N}}|\dot{u}|^{2}(t,x)dx$ and $\int^{t}_{0}\int_{\T^{N}}f{N}(s)|\n\dot{u}|^{2}dxds$. First we rewrite the momentum equation on the following form:
$$\rho\dot{u}-\mu\D u-(\lambda+\mu)\n{\rm div}u+\n P(\rho)=\rho f.$$
We apply to the momentum equation the operator $\frac{d}{dt}=\p_{t}+u\cdot\n$, we recall here the following equalities:
$$
\begin{aligned}
\frac{d}{dt}\rho\dot{u}^{j}&=\rho \frac{d}{dt}\dot{u}^{j}+\p_{t}\rho \dot{u}^{j}+\dot{u}^{j}\sum_{k}\p_{k}\rho \;u^{k},\\
&=\rho \frac{d}{dt}\dot{u}^{j}-\rho{\rm div}u \dot{u}^{j} ,
\end{aligned}
$$
We have next:
$$
\begin{aligned}
\mu\frac{d}{dt}\D u^{j}&=\mu\p_{t}\D u^{j}+\sum_{k}\p_{k}\D u^{j} u^{k},\\
&=\e(\p_{t}\D u^{j}+{\rm div}(\D u^{j} u)-\D u^{j}{\rm div}u),
\end{aligned}
$$
and (where $D={\rm div}u$):
$$
\begin{aligned}
(\lambda+\mu)\frac{d}{dt}\p_{j}{\rm div}u&=(\lambda+\mu)(\p_{t}\p_{j}D+{\rm div}(\p_{j}D u)-\p_{j}D{\rm div}u),
\end{aligned}
$$
We obtain finally:
\begin{equation}
\begin{aligned}
&\rho\frac{d}{dt}\dot{u}^{j}+\p_{j}\p_{t}P(\rho)+{\rm div}\p_{j}P(\rho)=\mu(\p_{t}\D u^{j}+{\rm div}(\D u^{j} u)\\
&\hspace{7cm}+
(\lambda+\mu)(\p_{t}\p_{j}D+{\rm div}(\p_{j}D u)).
\end{aligned}
\label{b2.10}
\end{equation}
We shall make use of the following transport theorem if $\rho\dot{w}=f_{1}$ and if $h=g(t)$, then:
$$\int^{t}_{0}\int_{\T^{N}} \frac{1}{2} \p_{s}(h\rho w^{2})dxds=\int^{t}_{0}\int_{\T^{N}} (\frac{1}{2}h^{'}\rho w^{2}+hwf)dxds.$$
We apply the previous result with:
$$f_{1}=-\p_{j}\p_{t}P(\rho)-{\rm div}\p_{j}P(\rho)+\mu(\p_{t}\D u^{j}+{\rm div}(\D u^{j} u)+
(\lambda+\mu)(\p_{t}\p_{j}D+{\rm div}(\p_{j}D u)),$$
and with $h(s)=f(s)^{2}$, we obtain then:
\begin{equation}
\begin{aligned}
&\frac{1}{2}f(t)^{2}\int\rho(t,x)|\dot{u}(t,x)|^{2}dx=\int^{t}_{0}\int_{\T^{N}}\frac{1}{2} f(s)
f^{'}(s)\rho|\dot{u}|^{2}
dxds+\int^{t}_{0}\int_{\T^{N}} f(s)^{2}\dot{u}^{j}\\
&\times[-(\p_{j}\p_{t}P+{\rm div}(\p_{j}Pu)+\mu[\D\p_{t}u^{j}+{\rm div}(\D u^{j}u)]
+(\lambda+\mu)[\p_{j}\p_{t}D\\
&\hspace{10cm}+{\rm div}(\p_{j}D u)]dxds.
\end{aligned}
\label{12.11}
\end{equation}
Since $f(s)
f^{'}(s)\leq f(s)$ we can apply (\ref{5imp}) to bound the first term on the right.
Next by integrations by part we get:
$$
\begin{aligned}
-\int_{0}^{t}\int_{\T^{N}} f(s)^{2}\dot{u}^{j}(\p_{j}\p_{t}P+{\rm div}(\p_{j}Pu))dsdx& =
\int_{0}^{t}\int_{\T^{N}} f(s)^{2}(\p_{j}\dot{u}^{j}\p_{t}P+\p_{k}\dot{u}^{j}\p_{j}P u^{k}dsdx,\\
&=\int\int f(s)^{N}P^{'}(\p_{j}\dot{u}^{j}\p_{t}\rho+\p_{k}\dot{u}^{j}\p_{j}\rho u^{k})dsdx
\end{aligned}
$$
$$
\begin{aligned}
&\int_{0}^{t}\int_{\T^{N}} f(s)^{2}P^{'}[-\p_{j}\dot{u}^{j}(\rho\p_{k}u^{k}+\p_{k}\rho u^{k})+\p_{k}\dot{u}^{j}\p_{j}\rho u^{k}]dxds,\\
&=-\int_{0}^{t}\int_{\T^{N}} f(s)^{2}[P^{'}\rho D\p_{j}\dot{u}^{j}+\p_{k}Pu^{k}\p_{j}\dot{u}^{j}-\p_{j}Pu^{k}\p_{k}\dot{u}^{j}]dxds,\\
&=-\int_{0}^{t}\int_{\T^{N}}  f(s)^{2}[P^{'}\rho D\p_{j}\dot{u}^{j}-P(D\p_{j}\dot{u}^{j}-\p_{j}u^{k}\p_{k}\dot{u}^{j})]dxds.
\end{aligned}
$$
This term is therefore bounded in absolute value by:
$$
\begin{aligned}
&C\big(\int^{t}_{0}\int_{\T^{N}}P(\rho)f(s)^{2}|\n u|^{2}dxds\big)^{\frac{1}{2}}\big(\int^{t}_{0}\int_{\T^{N}}f(s)^{2}|\n \dot{u}|^{2}dxds\big)^{\frac{1}{2}}\\
&\hspace{2cm}\leq C_{\e}\int^{t}_{0}\int_{\T^{N}} P(\rho)f(s)^{2}|\n u|^{2}(s,x)dxds+\e\int^{t}_{0}\int_{\T^{N}} f(s)^{N}|\n \dot{u}(s,x)|^{2}dxds.
\end{aligned}
$$
The third term on the right side of (\ref{12.11}) may be written:
$$
\begin{aligned}
&-\mu\int^{t}_{0}\int_{\T^{N}} f^{2}[\n\dot{u}^{j}\cdot\n u^{j}_{t}+(\n\dot{u}^{j}\cdot u)\D u^{j}]dxds\\
&\hspace{1cm}=-\mu\int^{t}_{0}\int_{\T^{N}} f^{2}[\n\dot{u}^{j}\cdot(\n u^{j}_{t}+\n(\n\dot{u}^{j}\cdot u))+\dot{u}^{j}_{k}(u^{k}u^{j}_{ll}-(u^{j}_{l}u^{l})_{k})]dxds,\\
&\hspace{1cm}=-\mu\int^{t}_{0}\int_{\T^{N}}f^{2}[\n\dot{u}^{j}\cdot(\n u^{j}_{t}+\n(\n\dot{u}^{j}\cdot u))+\dot{u}^{j}_{k}(u^{k}u^{j}_{ll}-u^{j}_{lk}u^{l}-u^{j}_{l}u^{l}_{k})]dxds,\\
&\hspace{1cm}\leq -\mu\int^{t}_{0}\int_{\T^{N}} f^{2}|\n\dot{u}^{j}|^{2}+M\int\int f^{2}|\n u|^{2}(|\n\dot{u}|+|\dot{D}|)dxds.
\end{aligned}
$$
The last term on the right side of (\ref{12.11}) may be bound as follows:
$$
\begin{aligned}
&-(\lambda+\mu)\int^{t}_{0}\int_{\T^{N}}f^{2}(\p_{j}\p_{t}D+{\rm div}(\p_{j}D u) dxds\leq
-(\lambda+\mu)\int^{t}_{0}\int_{\T^{N}}f^{2} \dot{D}^{2}dxds\\
&\hspace{7cm}+M\int^{t}_{0}\int_{\T^{N}} f^{2}|\n u|^{2}(|\n\dot{u}|+|\dot{D}|)dxds.
\end{aligned}
$$
It then follows by Young's inequalities that:
\begin{equation}
\begin{aligned}
&f(t)^{2}\int_{\T^{N}}\rho|\dot{u}|^{2}(t,x)dx+\int^{t}_{0}\int_{\T^{N}}f^{2}(s)(\mu|\n\dot{u}|^{2}+(\lambda+\mu)|\dot{D}|^{2})dxds\\
&\leq M\big[C_{0}+C_{\e}C_{0}\|P(\rho)\|_{L^{\infty}}+
\int^{t}_{0}\int_{\T^{N}}\frac{1}{2}f(s)
f^{'}(s)\rho|\dot{u}|^{2}dxds\\
&\hspace{7cm}+\int^{t}_{0}\int_{\T^{N}}f^{2}(s)|\n u|^{4}dxds],\\
&\leq M\big[C_{0}+C_{\e}C_{0}\|P(\rho)\|_{L^{\infty}}+
A(t)+\int^{t}_{0}\int_{\T^{N}}f^{2}(s)|\n u|^{4}dxds].
\end{aligned}
\label{aa1}
\end{equation}
Next from the momentum equation, we obtain as in the works of D. Hoff in \cite{}:
$$\mu|\n w|^{2}=\mu({\rm div}(w\n w)+\p_{j}(\rho w\dot{u}^{k})-\p_{k}(\rho w\dot{u}^{j})+\rho(\dot{u}^{j}\p_{k}w-\dot{u}^{k}\p_{j}w).$$
Integrating we thus obtain:
\begin{equation}
\begin{aligned}
\int^{t}_{0}\int_{\T^{N}}f(s)|\n w|^{2}dxds&\leq M(\int^{t}_{0}\int_{\T^{N}} f(s)\rho|\dot{u}|^{2}(s,x)dsdx+\int^{t}_{0}\int_{\T^{N}}f(s)\rho|w|^{2}(s,x)dsdx),\\
&\leq M(\int^{t}_{0}\int_{\T^{N}} f(s)\rho|\dot{u}|^{2}(s,x)dsdx+\|\rho\|_{L^{\infty}}),\\
&\leq M(A(t)+\|\rho\|_{L^{\infty}}).
\end{aligned}
\label{aa2}
\end{equation}
Similarly we obtain:
\begin{equation}
\begin{aligned}
&\sup_{0<t\leq T}f(t)^{2}\int |\n\omega|^{2}dx\leq M(\int_{\T^{N}} f^{2}(t)\rho|\dot{u}|^{2}(t,x)dsdx+\int_{\T^{N}}f(t)^{2}\rho|w|^{2}(t,x)dx)
\\
&\hspace{4cm}\leq M(\int_{\T^{N}} f^{2}(t)\rho|\dot{u}|^{2}(t,x)dsdx+C\|\rho\|_{L^{\infty}}),\\
&\leq M(C_{0}+C_{\e}C_{0}\|P(\rho)\|_{L^{\infty}}+
A(t)+\int^{t}_{0}\int_{\T^{N}}f^{2}(s)|\n u|^{4}dxds+C\|\rho\|_{L^{\infty}}).
\end{aligned}
\label{aa3}
\end{equation}
To complete the estimates (\ref{aa1}) and (\ref{aa3}), we will need of estimation on
$\int^{t}_{0}\int_{\T^{N}}f^{2}(s)|\n u|(s,x)^{4}dxdt$ and to conclude we will use (\ref{5imp}).
\subsubsection*{Control of $\int^{t}_{0}\int_{\T^{N}}f^{2}(s)|\n u|(s,x)^{4}dxdt$}
We have then:
\begin{equation}
\int^{t}_{0}\int_{\T^{N}} f^{2}|\n u|^{4}dsdx\leq \int^{t}_{0}\int_{\T^{N}} f^{2}\big((RG)^{4}+\omega^{4})(s,x)+f^{2}(s)RP(\rho)^{4}(s,x)\big)dxds.
\label{vorti}
\end{equation}
Let focus us on the case $N=3$. When $N=3$ we can apply Gagliardo-Nirenberg, let:
$$\|RG\|_{L^{4}}\leq C\|G\|_{L^{2}}^{\frac{5}{8}}\|\n G\|_{L^{6}}^{\frac{3}{8}}.$$
We have then by Young's inequality and the fact that $(\lambda+2\mu)\D G={\rm div}(\rho (\dot{u}+g))$:
\begin{equation}
\begin{aligned}
&\int^{t}_{0}\int_{\T^{N}}f^{2}(s)(RG)^{4}(s,x)dsdx\leq \int^{t}_{0}f^{2}(s)\|G\|^{\frac{5}{2}}_{L^{2}}\|\n G\|_{L^{6}}^{\frac{3}{2}}ds,\\
&\leq\int^{t}_{0}f^{2}(s)\|G\|_{L^{2}}^{\frac{5}{2}}(\|\n \dot{u}\|_{L^{2}}+\|f\|_{L^{6}})^{\frac{3}{2}}ds,\\
&\leq\e\int^{t}_{0}\int_{\T^{N}}f^{2}(s)|\n\dot{u}|^{2}dxds+C_{\e}\int^{t}_{0}f^{2}(s)\|G\|_{L^{2}}^{10}ds,\\
&\leq\e\int^{t}_{0}\int_{\T^{N}}f^{2}|\n\dot{u}|^{2}dxds+C_{\e}\int^{t}_{0}(f(s)\|G\|_{L^{2}})^{10}\frac{1}{f(s)^{7}}ds,\\
&\leq\e\int^{t}_{0}\int_{\T^{N}}f^{2}|\n\dot{u}|^{2}dxds+C_{\e}C (\sup_{0\leq s\leq t}f(s)\|G\|_{L^{2}})^{10}.
\end{aligned}
\label{a5imp}
\end{equation}
From (\ref{a5imp}) and (\ref{aa1}), (\ref{aa3}), we can conclude that:
\begin{equation}
\begin{aligned}
&B(t)f(t)^{2}\int_{\T^{N}}\rho|\dot{u}|^{2}(t,x)dx+\int^{t}_{0}\int_{\T^{N}}f^{2}(s)(\mu|\n\dot{u}|^{2}+(\lambda+\mu)|\dot{D}|^{2})dxds\\
&\hspace{3cm}\leq M\big[C_{0}+C_{\e}C_{0}\|P(\rho)\|_{L^{\infty}}+CC_{e}A(t)(1+A(t)^{9})].\\
\end{aligned}
\label{aa11}
\end{equation}
\begin{equation}
\begin{aligned}
&\sup_{0<t\leq T}f(t)^{2}\int |\n\omega|^{2}dx
\leq M(C_{0}+C_{\e}C_{0}\|P(\rho)\|_{L^{\infty}}+
A(t)C\|\rho\|_{L^{\infty}}\\
&\hspace{9cm}+\e B(t)+CC_{e}A(t)^{10}).
\end{aligned}
\label{aa33}
\end{equation}
We can remark that all the inequalities (\ref{5imp}), (\ref{aa2}), (\ref{aa11}) and (\ref{aa33}) depend on the control of $\|\rho\|_{L^{\infty}}$. In the following subsection, we want explain that we can control $\rho$ in $L^{\infty}$ in finite time.
\subsubsection*{Conclusion}
We will treat by simplicity only the case $N=3$.
We want here to explain how to control the norm $L^{\infty}$ of the density $\rho$ in finite time.
From (\ref{29}), we have:
$$
\begin{aligned}
\log(\rho(t,x))\leq&\log(\|\rho_{0}\|_{L^{\infty}})+C\|(\D)^{-1}{\rm div} m_{0}\|_{L^{\infty}}+C\|(\D)^{-1}{\rm div}(\rho u)\|_{L^{\infty}}\\
&+C \int^{t}_{0}\|[u_{j},R_{i}R_{j}](\rho u_{i})(s,\cdot)\|_{L^{\infty}}ds,
\end{aligned}
$$
From the previous section and the regularizing effects on $u$ we obtain:
$$\log(\rho(t,x))\leq C_{t}+C \exp(C\exp(\int^{t}_{0}\phi(\|\rho(s,\cdot)\|_{L^{\infty}})ds)).$$
Consequently by Gr\"onwall lemma, there exists $T_{0}>0$ such that for all $T<T_{0}$:
$$\|\rho\|_{L^{\infty}}\leq C.$$
\subsubsection*{Proof of (\ref{b1.22}) and (\ref{b1.23})}
We want get now solutions such that we control $\n u$ in $L^{1}(BMO)$. To do it, we need of additional regularity on the velocity.
We will use again the technics introduced by D. Hoff in \cite{Hoffn2}. The idea is to obtain some estimates by interpolation
in ``killing'' the coupling between pressure and velocity.
First , we mollify initial data satisfying the conditions of theorem \ref{theo1} and then appeal to the result of \cite{H3}
to otain a solution $(\rho,u)$ defined at least for small time. We want now derive estimates on the solution independent of the mollifier and depending only on the initial data.
\\
In the sequel we will treat only by simplicity the case $N=3$.
Fixing the local in time solution $(\rho,u)$ described above on the interval $[0,T_{0}]$ with $T_{0}>0$, we therefore assume throughout this section that $C^{-1}\leq \rho\leq C$ with $C>0$.
We define a differential operator ${\cal L}$ acting on functions $w:[0,T]\times\T^{N}\rightarrow\T^{N}$ by
$${\cal L}w=\p_{t}(\rho w)+{\rm div}(\rho u\otimes w)-\mu\D w-\lambda\n{\rm div}w,$$
and we define $w_{1}$ and $w_{2}$ by:
\begin{equation}
{\cal L}w_{1}=0,\;\;(w_{1})_{/t=0}=u_{0},\;\;\;{\cal L}w_{2}=-\n P(\rho),\;\;(w_{2})_{/t=0}=0.
 \label{g2.1}
\end{equation}
We observe here that by uniqueness $w_{1}+w_{2}=u$. Straightforward energy estimates then show that:
\begin{equation}
\sup_{0\leq t\leq T_{0}}\int_{\T^{N}}|w_{1}(t,x)|^{2}dx+\int^{T_{0}}_{0}\int_{\T^{N}}|\n w_{1}|^{2}dxdt\leq C\int_{\T^{N}}|u_{0}|^{2}dx
 \label{g2.2}
\end{equation}
and:
\begin{equation}
\sup_{0\leq t\leq T_{0}}\int_{\T^{N}}|w_{2}(t,x)|^{2}dx+\int^{T_{0}}_{0}\int_{\T^{N}}|\n w_{2}|^{2}dxdt\leq CT \sup_{0\leq t}|P(\rho(t,\cdot)|^{2}.
 \label{g2.3}
\end{equation}
We shall derive (\ref{b1.22}) and (\ref{b1.23}) simultaneously as consequences of estimates for the following quantities:
$$\sup_{0\leq t\leq T_{0}}t^{1-k}\int_{\T^{N}}|\n w_{1}(t,x)|^{2}dx+\int^{T_{0}}_{0}\int_{\T^{N}} t^{1-k}|\dot{w}_{1}|^{2}dxdt,$$
for $k=0,1$ and,
$$\sup_{0\leq t\leq1}\int|\n w_{2}(t,x)|^{2}dx+\int^{1}_{0}\int|\dot{w}_{2}|^{2}dxdt.$$
To derive these bounds, we multiply equation (\ref{g2.1}) for $w_{1}$ and $w_{2}$ by $\dot{w}_{1}$ and $\dot{w}_{2}$, respectively and integrate. The details which are nearly identical to those in the previous section are quite straightforward, the essential point being that the spatial gradient of $\rho$ must be avoided. Indeed the procedure of D. Hoff in \cite{Hoffn2} allows to ``kill'' the coupling between the velocity and the pressure. The results are that with $C^{'}>0$:
\begin{equation}
\begin{aligned}
 &\frac{1}{2}\mu(\tau^{k}\int|\n w_{1}(\tau,x)|^{2}dx\big|)^{\tau=t}_{\tau=0}+\int^{t}_{0}\int_{\T^{N}}\tau^{k}\rho|\dot{w}_{1}|^{2}dxd\tau\leq\\
&\hspace{2cm}\frac{1}{2}\mu k\int^{t}_{0}\int_{\T^{N}}\tau^{k-1}|\n w_{1}|^{2}dxd\tau+\int^{t}_{0}\int_{\T^{N}}\tau^{k}|\n w_{1}|^{2}|\n u|dxd\tau,
\end{aligned}
 \label{g2.4}
\end{equation}
and
\begin{equation}
\begin{aligned}
 &\frac{1}{2}\mu\int_{\T^{N}}|\n w_{2}(\tau,x)|^{2}dx+\int^{t}_{0}\int_{\T^{N}}\rho|\dot{w}_{2}|^{2}dxd\tau\leq\\
&\hspace{0,5cm}\int_{\T^{N}}P(\rho(t,\cdot)){\rm div}w_{2}(\cdot,x)dx\big|^{t}_{0}+\int^{t}_{0}\int_{\T^{N}}(|\n w_{2}|^{2}|\n u|+
|\n w_{2}||\n u|)dxd\tau ,
\end{aligned}
 \label{g2.5}
\end{equation}
By proceding exactly as in the previous section, we obtain the following results:
\begin{equation}
\sup_{0\leq t\leq T_{0}}\int_{\T^{N}}|\n w_{1}(t,x)|^{2}dx+\int^{T_{0}}_{0}\int_{\T^{N}}|\dot{w}_{1}|^{2}dxdt\leq C\|u_{0}\|_{H^{1}}^{2},
 \label{g2.7}
\end{equation}
\begin{equation}
\sup_{0\leq t\leq T_{0}}t\int_{\T^{N}}|\n w_{1}(t,x)|^{2}dx+\int^{T_{0}}_{0}\int t|\dot{w}_{1}|^{2}dxdt\leq C\|u_{0}\|_{L^{2}}^{2},
 \label{g2.8}
\end{equation}
\begin{equation}
\sup_{0\leq t\leq T_{0}}\int_{\T^{N}}|\n w_{2}(t,x)|^{2}dx+\int^{T_{0}}_{0}\int_{\T^{N}}|\dot{w}_{2}|^{2}dxdt\leq C C_{0},
 \label{g2.9}
\end{equation}
where $C_{0}$ depends only of the initial data. Now since the solution operator $u_{0}\longrightarrow w_{1}(t,\cdot)$ is linear, we can apply a standard Riesz-Thorin interpolation argument to deduce from (\ref{g2.7}) and (\ref{g2.8}) that:
\begin{equation}
\sup_{0\leq t\leq T_{0}}t^{1-\beta}\int^{\T^{N}}|\n w_{1}(t,x)|^{2}dx+\int^{T_{0}}_{0}\int_{\T^{N}} t^{1-\beta}|\dot{w}_{1}|^{2}dxdt\leq C\|u_{0}\|_{H^{\beta}}^{2}.
 \label{g2.10}
\end{equation}
As $u=w_{1}+w_{2}$, we then conclude from (\ref{g2.9}) and (\ref{g2.10}) that:
\begin{equation}
\sup_{0\leq t\leq T_{0}}t^{1-\beta}\int_{\T^{N}}|\n u(t,x)|^{2}dx+\int^{T_{0}}_{0}\int_{\T^{N}} t^{1-\beta}|\dot{u}|^{2}dxdt\leq C
C_{0}\|u_{0}\|_{H^{\beta}}^{2}.
 \label{g2.11}
\end{equation}
The next step is to obtain bounds for the terms
$$\sup_{0\leq t\leq T_{0}}t^{2-\beta}\int_{\T^{N}}|\n u(t,x)|^{2}dx+\int^{T_{0}}_{0}\int_{\T^{N}} t^{2-\beta}|\n \dot{u}|^{2}dxdt$$
appearing in (\ref{b1.22}). To do this, we multiply the momentum equation of (\ref{1})  by $t^{2-\beta}\dot{u}$ and integrate.
The details are exactly as in the previous section, except now we apply the $\beta$ dependent smoothing rates established in (\ref{g2.10}). Combining
these bounds with (\ref{g2.2}), (\ref{g2.3}) and (\ref{g2.10}), we then obtain (\ref{b1.22}) for times $t\leq T_{0}$. 
To prove (\ref{b1.23}), we observe that for $k=0,1$,
$$\sup_{0\leq t\leq 1}\|w_{1}(t,\cdot)\|_{H^{k}}\leq C\|u_{0}\|_{H^{k}},$$
by (\ref{g2.2}) and (\ref{g2.7}). Thus:
$$\sup_{0\leq t\leq 1}\|w_{1}(t,\cdot)\|_{H^{\beta}}\leq C\|u_{0}\|_{H^{k}},$$
for $\beta\in [0,1]$. As $u=w_{1}+w_{2}$, and applying (\ref{g2.9}) we obtain that:
$$\sup_{0\leq t\leq 1}\|w_{1}(t,\cdot)\|_{H^{\beta}}\leq C C_{0},$$
and then for $r\in(2,\frac{6}{3-2\beta})$ in the case that $\beta>0$, that:
$$\sup_{0\leq t\leq 1}\|u(r,\cdot)-\widetilde{u}\|_{L^{r}}\leq C C_{0}.$$
This proves (\ref{b1.23}).
\subsubsection*{Regularity on the gradient of the velocity $u$ in $L^{1}_{T_{0}}(W^{1,\alpha}+BMO)$ with $\alpha>N$}
\label{v1reg}
Here we want examine the regularity of the gradient of the velocity and prove that $\n u$ is in $L^{1}_{T_{0}}(BMO)$ to prove (\ref{c1.23}). More precisely we will verify that the new variable $v_{1}$ introduced in \cite{H3} called ``effective velocity' belongs in $L^{1}_{T_{0}}(W^{2,\alpha})$ with $\alpha>N$, let $\n v_{1}\in L^{1}{T_{0}}(L^{\infty})$. 
We recall here the definition of $v_{1}$ introduced in \cite{H3}.
The idea of \cite{H3} was to introduce a variable $v_{1}$ which allows to kill the coupling between the velocity and the pressure in the momentum equation of (\ref{1}).
In this goal, we need to integrate the pressure term in the study of the linearized equation of the momentum equation. To do this, we will try to express the gradient of the pressure as a Laplacian term, so we set:
$$\D v=\n P(\rho).$$
We have then $v=(\D)^{-1}\n P(\rho)$ with $(\D)^{-1}$ the inverse Laplacian with zero value on $\T^{N}$. In the sequel we will set:
$$v_{1}=u-\frac{\lambda+2\mu}v.$$
We have then:
\begin{equation}
 \begin{aligned}
\D u=&\n F+{\rm div}\omega+(2\mu+\lambda)^{-1}\n (P(\rho)),\\
&=\D v_{1}+(2\mu+\lambda)^{-1}\n (P(\rho)).
\end{aligned}
\label{11.22}
\end{equation}
We can easily show that $\int^{T_{0}}_{0}\|\n v\|_{L^{\infty}}dt<+\infty$ if (\ref{b1.22}) holds.
To see this, we apply standard elliptic theory 
on $v_{1}$ and the fact that $\D G={\rm div}(\rho\dot{u}-\rho g)$. We will use in particular the fact that ${\rm div}v_{1}=(\lambda+2\mu)G$ and ${\rm curl}v_{1}=\omega$.
We consider here only the case $N=3$. The case $N=2$ follows the same lines. For some $\alpha>3$ and $\e>0$ determined by $\alpha$ we have then:
$$
\begin{aligned}
\|\n v_{1}\|_{L^{\infty}}&\leq C(\|\n F\|_{L^{\alpha}}+\|\n\omega\|_{L^{\alpha}}), \\
&\leq C(\|\rho\dot{u}(t,\cdot)\|_{L^{\alpha}}+\|\rho g\|_{L^{\alpha}}+\|\n\omega\|_{L^{\alpha}}),\\
&\leq C(\|\rho\|_{L^{\infty}}\|\dot{u}(t,\cdot)-\widetilde{\dot{u}}\|^{\frac{1-\e}{2}}_{L^{2}}\|\n\dot{u}(t,\cdot)\|^{\frac{1+\e}{2}}_{L^{2}}+
\|\rho\|_{L^{\infty}}\|g\|_{L^{\alpha}}+\|\n\omega\|_{L^{\alpha}}),\\
&\leq C(\|\rho\|_{L^{\infty}}(\|\dot{u}(t,\cdot)\|^{\frac{1-\e}{2}}_{L^{2}}+\widetilde{\dot{u}}^{\frac{1-\e}{2}})\|\n\dot{u}(t,\cdot)\|^{\frac{1+\e}{2}}_{L^{2}}+
\|\rho\|_{L^{\infty}}\|g\|_{L^{\alpha}}+\|\n\omega\|_{L^{\alpha}}),
\end{aligned}
$$
We recall here that $\widetilde{\rho\dot{u}}=\widetilde{\rho g}$, we have then as $\frac{1}{\rho}\geq C$ on $[0,T_{0}]$:
$$\widetilde{u}\leq\frac{1}{\|\rho\|_{L^{\infty}_{T_{0}}}}\widetilde{\rho g}.$$
so that by using (\ref{b1.22}), we obtain:
$$
\begin{aligned}
\int^{T_{0}}_{0}\|\n v_{1}\|_{L^{\infty}}dt&\leq C\int^{T_{0}}_{0}t^{\beta}(t^{1-s}\int_{\T^{N}}|\dot{u}|^{2}dx)^{\frac{1-\e}{4}}(t^{\sigma}\int_{\T^{N}}|\n\dot{u}|^{2}dx)^{\frac{1+\e}{4}}dt+C_{0},\\
&\leq C(\int^{T_{0}}_{0}t^{2\beta}dt)^{\frac{1}{2}}+C_{0},
\end{aligned}
$$
 with $s=\N+\e-1$ ($\e>0$) and where $4\beta=(s-1)(1-\e)-(\sigma+\e)$. The above integral is therefore finite as $2\beta>-1$. 
A similar result result holds for $N=2$ with $s>0$. Thus for the solution constructed in the previous section, $\int^{T}_{0}\|\n v_{1}(t,\cdot)\|_{L^{\infty}}dt$ is finite if (\ref{b1.22}) holds, $\inf\rho_{0}\geq c>0$ and $u_{0}\in H^{\N+\e-1}$,with $\e>0$.\\
More precisely we have proved in fact that:
\begin{equation}
 \n v_{1}\in L^{1}_{T}(W^{1,\alpha})\hookrightarrow L^{1}_{T}(B^{1+\e}_{N,1}).
\label{gainv1}
\end{equation}
with $\alpha=N+2\e$ where $\e>0$.\\
As $P(\rho)\in L^{\infty}$ we deduce from (\ref{11.22}) and the results of Calderon-Zygmund, that:
$$\n u \in L^{1}_{T_{0}}(BMO).$$
\section{Proof of theorem \ref{theo1}}
\label{section5}
The above arguments are not rigorous, since we have to assume that $\rho(t,x)$ is positive for all $(t,x)$. In order to deal with possibly vanishing densities, we remark that if we assume as in \cite{H3} that $\rho_{0}$ is also bounded from below, we can get $L^{\infty}$ bounds for $\log\rho$. In that case, when $N=2,3$, there exists $T_{0}>0$ such that for all $t<T_{0}$:
\begin{equation}
\|\rho\|_{L^{\infty}((0,t)\times\T^{N})}+\|\frac{1}{\rho}\|_{L^{\infty}((0,t)\times\T^{N})}\leq C_{t}.
 \label{119}
\end{equation}
Thus $\rho$ is also bounded from below for small enough times, so that vacuum does not form on $[0,T_{0}]$. It follows that starting from a general bounded initial density $\rho_{0}$, we can apply the above arguments to a weak solution $(\rho^{n},u^{n})$ corresponding to initial values $\rho_{0}^{n}=\rho_{0}+\frac{1}{n}$ and $u_{0}^{n}=u_{0}$ which converge strongly in $L^{\infty}(\T^{N})$ to $\rho_{0}$ and $u_{0}$. In view of the weak stability results of system (\ref{1}) given by E. Feireisl et al in \cite{5F3}, $\rho^{n}$ and $\sqrt{\rho^{n}}u^{n}$ respectively converge to $\rho$ and $u$ in $C([0,T_{0}],L^{q}(\T^{N}))$ and $L^{2}((0,T_{0})\times\T^{N})$ for all $q<\gamma-1+\frac{2\gamma}{N}$. Hence the uniform $L^{\infty}$ bounds on $\rho^{n}$ yield $L^{\infty}$ bounds on $\rho$.\\
In the above formal derivations, we assumed that there exists global weak solutions of (\ref{1}). This problem does
not occur when we take $\gamma$ larger than $\frac{N}{2}$ via the works of Feireisl et al in \cite{5F3}. When $1<\gamma\leq\frac{N}{2}$, we can approximate solutions of (\ref{1}) by a global weak solutions of (\ref{1}) by a global weak solution $(\rho^{n},u^{n})$ corresponding to a modified pressure law that satisfies:
\begin{equation}
P_{n}(\rho)=P(\rho)+\frac{1}{n}\rho^{2},
 \label{120}
\end{equation}
and the same initial data $(\rho_{0},u_{0})$. Applying the above arguments on $(\rho^{n},u^{n})$, we obtain all the uniform bounds for $\rho^{n}$ and $u^{n}$ on $[0,T_{0}]$, where $T_{0}$ does not depend on $n$.
As a result, we also have uniform $L^{2}((0,T_{0})\times\T^{N})$ bounds on $\n u^{n}$ and $L^{\infty}(0,T_{0},L^{2}(\T^{N}))$ bounds on $\sqrt{\rho}^{n}u^{n}$. Hence the weak stability results hold since the initial data are not functions of $n$ and $(\rho^{n},u^{n})$ converge to $(\rho,u)$ in $L^{2}((0,T_{0})\times\T^{N})^{N+1}$, where $(\rho,u)$ is a solution of (\ref{1}) on $(0,T_{0})$. We refer to \cite{13} for complete details of the stability proof. Let us emphasize that one of the key arguments for proving weak stability of solutions of (\ref{1}) is to obtain uniform $L^{2}((0,T_{0})\times\T^{N})$ bounds on $\rho^{n}$ to renormalize the transport equation. This is the case here as we have uniforms bounds on the density in $L^{\infty}((0,T_{0})\times\T^{N})$.
\section{Proof of corollary \ref{corollaire1}}
\label{section6}
\subsubsection{Control of $\rho\in L^{\infty}(B^{\e}_{\infty,\infty})$}
In view of proposition \ref{transport} where in our case $h(\rho)=P(\rho)$, we have for all $t\in[0,T]$ and $0<\e<1$:
\begin{equation}
 \|\rho\|_{\widetilde{L}^{\infty}_{t}(B^{\e}_{\infty,\infty})}\leq e^{CV(t)}(1+\|\rho_{0}\|_{B^{\e}_{\infty,\infty}}),
\label{38}
\end{equation}
where $V(t)=\int^{t}_{0}\big(\|\n u(\tau)\|_{L^{\infty}}+\|{\rm div}v_{1}(\tau)\|_{B^{\e}_{\infty,\infty}}+\|\rho(\tau)\|^{s}_{L^{\infty}}\big)d\tau$,
where $s$ the smallest integer such that $P^{'}\in W^{s,\infty}$. We have seen by (\ref{gainv1}) that $\n v_{1}\in L^{1}(0,T, B^{\e}_{\infty,\infty})$ with $\e>0$, the main difficulty is to control $\n u\in L^{1}(0,T,L^{\infty})$, for this we recall that:
$$
\begin{aligned}
\|\n u\|_{L^{1}_{T}(L^{\infty})}&\leq \|\n v_{1}\|_{L^{1}_{T}(B^{\e}_{\infty,\infty})}+\|P(\rho)\|_{L^{1}_{T}(B^{0}_{\infty,1})},\\
&\leq \|\n v_{1}\|_{L^{1}_{T}(B^{\e}_{\infty,\infty})}+\|\rho\|^{s}_{L^{\infty}_{T}(L^{\infty})}\|\rho\|_{L^{1}_{T}(B^{0}_{\infty,1})}.
\end{aligned}
$$
Unsurprisingly the result comes from the well known following estimate:
Next we have:
$$
\begin{aligned}
\|\rho(t)\|_{B^{0}_{\infty,1}}\leq C\|\rho(t)\|_{B^{0}_{\infty,\infty}}\log (e+\frac{\|\rho(t)\|_{B^{\e}_{\infty,\infty}}}{\|\rho(t)\|_{B^{0}_{\infty,\infty}}}),
\end{aligned}
$$
and we recall the following inequality:
$$\forall x>0, \;\forall\delta>0,\;\log(e+\frac{\delta}{x})\leq\log(e+\frac{1}{x})(1+\log\delta).$$
We obtain then from the previous inequality
$$\|\rho(t)\|_{B^{0}_{\infty,1}}\leq \|\rho(t)\|_{B^{0}_{\infty,\infty}}\big(1+\log(\|\rho(t)\|_{B^{\e}_{\infty,\infty}})\big)\log (e+\frac{1}{\|\rho(t)\|_{B^{0}_{\infty,\infty}}}),$$
Let $X(t)=\int^{t}_{0}\|\rho(s)\|_{B^{0}_{\infty,1}}ds$, we have then:
$$V(t)\leq C\big( X(t)+\int^{t}_{0}\big(\|\n v_{1}(\tau)\|_{L^{\infty}}+\|{\rm div}v_{1}(\tau)\|_{B^{\e}_{\infty,\infty}})d\tau\big).$$
Combining (\ref{38}) and the previous inequality leads to:
$$
\begin{aligned}
X(t)&\leq \int^{t}_{0}\|\rho(s)\|_{B^{0}_{\infty,\infty}}\big(1+CV(t)+\log(1+\|\rho_{0}\|_{B^{\e}_{\infty,\infty}})\big)\log (e+\frac{1}{\|\rho(s)\|_{B^{0}_{\infty,\infty}}})ds,\\
&\leq \int^{t}_{0}\|\rho(s)\|_{B^{0}_{\infty,\infty}}\big(1+CX(t)+C\int^{t}_{0}\big(\|\n v_{1}(\tau)\|_{L^{\infty}}+\|{\rm div}v_{1}(\tau)\|_{B^{\e}_{\infty,\infty}}\big)d\tau\\
&\hspace{5cm}+\log(1+\|\rho_{0}\|_{B^{\e}_{\infty,\infty}})\big)\log (e+\frac{1}{\|\rho(s)\|_{B^{0}_{\infty,\infty}}})ds.
\end{aligned}
$$
Applying Gronwall inequality and inequality (\ref{gainv1}) shows that:
$$
\begin{aligned}
X(t)&\leq C_{t,0}\exp(C\int^{t}_{0}\|\rho(s)\|_{B^{0}_{\infty,\infty}} \log (e+\frac{1}{\|\rho(s)\|_{B^{0}_{\infty,\infty}}})ds),\\
&\leq C_{t,0}\exp(C\int^{t}_{0}(1+\|\rho(s)\|_{L^{\infty}})ds),
\end{aligned}
$$
where $C_{t,0}$ depnds only of the time $t$ and the initial data. As $\rho\in L^{\infty}_{t}(L^{\infty})$, we conclude that $X(t)\leq C_{t}$ and by this way we have proved that:
\begin{equation}
 \|\rho\|_{L^{\infty}_{T}(B^{\e}_{\infty,\infty})}\leq C_{0,T},
\label{gainclas}
\end{equation}
where $C_{t,0}$ depnds only of the time $T$ and the initial data.
\subsubsection{Control of $\rho\in L^{1}(B^{1}_{N,1})$}
\label{section61}
In this case, we need to show for the sequel that $\rho\in L^{\infty}(B^{1}_{N,1})$, and for this we proceed exactly as previous.
Indeed as $\rho_{0}\in B^{1}_{N,1}$, by proceding as in the previous section we can show that $\rho\in L^{\infty}_{T}(B^{1}_{N,1})$. Next by using proposition \ref{transport}, we obtain the fact that:
\begin{equation}
\|\rho\|_{L_{T}^{\infty}(B^{1}_{N,1})}\leq e^{CV(T)}(1+\|q_{0}\|_{B^{1}_{N,1}}),
\label{gainclass1}
\end{equation}
with $V(T)=\int^{T}_{0}(\|\n v_{1}(\tau)\|_{B^{1+\e}_{\infty,\infty}}+\|q(\tau)\|_{B^{\e}_{\infty,\infty}})d\tau.$
\subsection{Uniqueness}
We now discuss the uniqueness of the solutions of theorem \ref{theo1}.
For this we want use the result of P. Germain \cite{PG} which is a result of weak-strong uniqueness. In the sequel we will note $(\rho_{1},u_{1})$ the solution of the theorem \ref{theo1} which exits on the time interval $[0,T_{0}]$.
We have shown that our solution check $\rho\in L^{\infty}(L^{\infty})$. By theorem \ref{theo1}, we obtain that our solution verify the following inequalities:
\begin{equation}
 \begin{aligned}
&\sup_{0<t\leq+\infty}\int_{\T^{N}}[\frac{1}{2}\rho(t,x)|u(t,x)|^{2}+|P(\rho(t,x))|+\sigma(t)|\n u(t,x)|^{2}dx\\
&+\sup_{0<t\leq +\infty}\int_{\T^{N}}[\frac{1}{2}\rho(t,x)f(t)^{N}(\rho|\dot{u}(t,x)|^{2}+|\n \omega(t,x)|^{2})dx\\
&+\int^{+\infty}_{0}\int_{\T^{N}}[|\n u|^{2}+f(s)\rho|\dot{u}|^{2}+|\omega|^{2})+\sigma^{N}|\n \dot{u}|^{2}]dxdt\\
&\hspace{9cm}\leq C(C_{0}+C_{f})^{\theta},
 \end{aligned}
\label{51.21}
\end{equation}
and we obtain moreover:
\begin{equation}
\begin{cases}
 \begin{aligned}
&\sqrt{\rho}\p_{t}u\in L^{2}_{t}(L^{2}(\T^{N})),\\
&\sqrt{t}{\cal P}u\in L^{2}_{T}(H^{2}(\T^{N})),\\
&\sqrt{t}G=\sqrt{t}[(\lambda+2\mu){\rm div}u-P(\rho)]\in L^{2}_{T}(H^{1}(\T^{N})),\\
&\sqrt{t}\n u\in L^{\infty}_{T}(L^{2}(\T^{N})),
 \end{aligned}
\end{cases}
\label{513}
\end{equation}
Now by the result of P. Germain in \cite{PG}, we are able to prove that $(\rho,u)=(\rho_{1},u_{1})$ on $[0,T_{0}]$. To see this we have just to verify that $(\rho_{1},u_{1}$ verify the conditions of the theorem 2.2 of \cite{PG}.
For simplicity we prove only the result for $N=3$. We know that $\n\rho_{1}\in L^{\infty}(B^{\e}_{N,1})\hookrightarrow L_{T_{0}}^{\infty}(L^{N})$ and $\n u_{1}\in L^{1}_{T_{0}}(L^{\infty})$.
The main thing is to see that $\sqrt{t} \dot{u}_{1}\in L^{2}_{T_{0}}(L^{N})$.
\\
We recall by Gagliardo-Nirenberg inequalities that:
$$\sqrt{t}\|\dot{u}_{1}\|_{L^{3}}\leq(t^{\frac{1}{4}-\frac{\e}{2}}\|\dot{u}_{1}\|_{L^{2}})^{\frac{1}{2}}(t^{\frac{3}{2}-\frac{\e}{2}}
\|\n \dot{u}_{1}\|_{L^{2}})^{\frac{1}{2}}t^{\frac{\e}{2}}.$$
From the inequalities (\ref{b1.22}), we deduce that $\sqrt{t}\dot{u}_{1}\in L^{2}(L^{3})$.
\section{Proof of theorem \ref{theo3}}
\label{section7}
\subsection{How to obtain a regularizing effect on $v_{1}$ when $\rho\in L^{\infty}(L^{q})$}
For the moment, we don't give conditions on $q$, we will precise his value in the sequel of the proof.\\
We want now use our change of variable $v_{1}$ introduced in the previous sections. The interest of this new variable is to ''kill'' in a certain way the coupling velocity-pressure. In this goal, we can now rewrite the momentum equation of system (\ref{1}). We obtain then the following equation where we have set $\nu=2\mu+\lambda$:
$$\rho\p_{t}u+\rho u\cdot \n u-\mu\D\big(u-\frac{1}{\nu}v\big)-(\lambda+\mu)\n{\rm div}\big(u-\frac{1}{\nu}v\big)=\rho g,$$
where we recall that $v=(\D)^{-1}(\n P(\rho))$ with $(\D)^{-1}$ the inverse Laplacian with zero mean value on $\T^{N}$.
As $v_{1}=u-\frac{1}{\nu}v$  we have:
$$\rho\p_{t}v_{1}+\rho u\cdot \n u-\mu\D v_{1}-(\lambda+\mu)\n{\rm div}v_{1}=\rho g-\frac{1}{\nu}\rho\p_{t}v.$$
As ${\rm div}v=P(\rho)-\int_{\T^{N}}P(\rho)dx$, from the transport equation we obtain:
$$
\begin{aligned}
{\rm div}\p_{t}v&=-P^{'}(\rho)\rho{\rm div}u-\n P(\rho)\cdot u+\widetilde{P^{'}(\rho)\rho{\rm div}u}+\widetilde{\n P(\rho)\cdot u}\\
&=-{\rm div}(P(\rho)u)+(P(\rho)-\rho P^{'}(\rho)){\rm div}u-\widetilde{P(\rho){\rm div}u}+\widetilde{\rho P^{'}(\rho)){\rm div}u}.
\end{aligned}
$$
In the sequel we will need to use the Bogovskii operator that we note $\Lambda^{-1}$ (see \cite{NS} p168 for a definition), we obtain then:
\begin{equation}
 \p_{t}v=\Lambda^{-1}\big(-{\rm div}(P(\rho)u)+(P(\rho)-\rho P^{'}(\rho)){\rm div}u-\widetilde{P(\rho){\rm div}u}+\widetilde{\rho P^{'}(\rho)){\rm div}u}\big).
\label{Bogo}
\end{equation}
We get finally
\begin{equation}
\rho\p_{t}v_{1}-\mu\D v_{1}-(\lambda+\mu)\n{\rm div}v_{1}=\rho g-\rho u\cdot \n u-\frac{1}{\nu}\rho\p_{t}v.
 \label{hchaleur}
\end{equation}
We set $f(t)=\min(t,1)$ and we remark that $f(0)=0$. We multiply then (\ref{hchaleur}) by $f(t)\p_{t}v_{1}$ and integrating on $(0,T)\times\T^{N}$ we obtain where $\xi=\mu+\lambda$:
\begin{equation}
\begin{aligned}
&\int^{t}_{0}\int_{\T^{N}}f(s)\rho|\p_{s}v_{1}|^{2}dxds+\frac{1}{2}\int_{\T^{N}}f(t)\big(\mu|\n v_{1}(t,x)|^{2}+\xi({\rm div} v_{1})^{2}(t,x)\big)dx\leq\\
&\hspace{3cm}\int^{t}_{0}\int_{\T^{N}}f^{'}(s)(\mu|\n v_{1}(t,x)|^{2}+\xi({\rm div} v_{1})^{2}(t,x)\big)dxds\\
&\hspace{1cm}+\int^{t}_{0}\int_{\T^{N}}
\rho u\cdot\n u \, f(s)\p_{t}v_{1}dxds+\int^{t}_{0}\int_{\T^{N}}(\rho g-\frac{1}{\nu}\rho\p_{s}v)f(s)\p_{s}v_{1}dxds.
\end{aligned}
 \label{h74}
\end{equation}
We have next to control the terms on the right hand side of (\ref{h74}). In the sequel we will for simplicity treat only the case $N=3$.
We can recall now that from the works of Mellet and Vasseur in \cite{5MV2}, we control the velocity $u$ in $L^{\infty}(L^{\infty})$
as we have assume that $\rho\in L^{\infty}(L^{3\gamma+\e})$ with $\e>0$. In fact from the inequality (\ref{againvitesse}), we could get a gain on $\rho^{\frac{1}{p}}u\in L^{\infty}(L^{p})$ with $p$ arbitraly big if $P(\rho)\in L^{p}(L^{\frac{3p}{p+1}})$ and in this case
$\rho\in L^{\gamma p}(L^{\frac{3\gamma p}{p+1}})$. We begin now with:
$$
\begin{aligned}
&\int^{t}_{0}\int_{\T^{N}}
u\cdot\n u \,\rho f(s)\p_{s}v_{1}dxds=\int^{t}_{0}\int_{\T^{N}}
\big(u\cdot\n v_{1}+u\cdot\n(\D)^{-1}\n(P(\rho)) \big)\,\rho f(s)\p_{s}v_{1}dxds.
\end{aligned}
$$
\subsection*{Estimates on the term $\int^{t}_{0}\int_{\T^{N}}u\cdot\n(\D)^{-1}\n(P(\rho))\,\rho f(s)\p_{t}v_{1}dxds$}
We start with treating the easiler term:
$$
\begin{aligned}
&\big|\int^{t}_{0}\int_{\T^{N}}u\cdot\n(\D)^{-1}\n(P(\rho))\,\rho f(s)\p_{t}v_{1}dxds\big|\leq\\ &C_{\e}\int^{t}_{0}\int_{\T^{N}}f(s)\rho |u\cdot\n(\D)^{-1}\n(P(\rho))|^{2}dxds+\e\|\sqrt{f(t)\rho}\p_{t}v_{1}\|^{2}_{L^{2}_{t}(L^{2}(\T^{N}))}
\end{aligned}
$$
As $\rho\in L^{\infty}(L^{q})$ and so $P(\rho)\in L^{\infty}(L^{\frac{q}{\gamma}})$, we have $\rho^{\frac{3\gamma-q}{q}}u\in L^{\infty}(L^{\frac{q}{3\gamma-q}})$ (we recall that here $q<3\gamma$). We have then to treat
$\sqrt{f(s)}\rho^{\frac{3q-6\gamma}{2q}}\rho^{\frac{3\gamma-q}{q}}|u| |\n(\D)^{-1}\n(P(\rho))|$ in $L^{2}_{t}(L^{2}(\T^{N})$.
We must have in this case by H\"older's inequalities:
\begin{equation}
\frac{3q-6\gamma}{2q^{2}}+\frac{\gamma}{q}+\frac{3\gamma-q}{q}=\frac{3q-6\gamma+8\gamma q-2q^{2}}{2q^{2}}\leq\frac{1}{2}\Longrightarrow 6q^{2}-2q(3+8\gamma)+12\gamma\geq 0.
\label{egal1}
\end{equation}
In this case, we get:
$$
\begin{aligned}
&\int^{t}_{0}\int_{\T^{N}}f(s)\rho |u\cdot\n(\D)^{-1}\n(P(\rho))|^{2}dxds\leq
\int^{t}_{0}f(s)\|\rho\|_{L^{\infty}(L^{q})}^{\frac{2q}{3\gamma-q}}\|\rho\|_{L^{\infty}(L^{q})}^{\frac{3q-6\gamma}{2q}}
\|\rho\|_{L^{\infty}(L^{q})}^{\gamma}ds\\
\end{aligned}
$$
And we can conclude.
\subsubsection*{Regularizing effect on $\D v_{1}$}
We want here using the regularizing effect on $v_{1}$. To do it we use the momentum equation (\ref{hchaleur}) and we have:
\begin{equation}
\begin{aligned}
\mu\D v_{1}+(\lambda+\mu)\n{\rm div}v_{1}=\rho\p_{t}v_{1}+\rho u\cdot \n v_{1}+\rho u\cdot\n(\D)^{-1}\n( P(\rho))+\rho g-\frac{\rho}{\nu}\p_{t}v.
\end{aligned}
 \label{bchaleur1}
\end{equation}
%
We want use the ellipticity of (\ref{bchaleur1}) to deduce regularizing effects on $v_{1}$.
\subsubsection*{Estimate on the term $\int^{t}_{0}\|\sqrt{f(s)\rho}u\cdot\n v_{1}\|_{L^{2}}^{2}ds$}
$$
\begin{aligned}
&\int^{t}_{0}\|\sqrt{f(s)\rho}u\cdot\n v_{1}\|_{L^{2}}^{2}ds,\\
&\leq\int^{t}_{0} \|\rho^{\frac{3\gamma-q}{q}}u\|_{L^{\frac{q}{3\gamma-q}}}^{2}\|\sqrt{\rho f(s)}\n v_{1}\|^{2}_{L^{q_{1}}}\|\rho\|_{L^{q}}^{\frac{3q-6\gamma}{q}}ds,\\
\end{aligned}
$$
Here by Gagliardo-Nirenberg, we get:
$$\|\sqrt{f(s)}\n v_{1}\|^{2}_{L^{q_{1}}}\leq \|\sqrt{f(s)}\D v_{1}\|^{2(1-\frac{\e}{2})}_{L^{p}}f(s)^{\frac{\e}{2}}\|\n v_{1}\|^{2\frac{\e}{2}}_{L^{2}},$$
with $\frac{1}{q_{1}}=\frac{1}{2}+\frac{1}{2q}-\frac{1}{N}-\frac{\e}{4q}-\frac{\e}{2N}$ and $q>\frac{N}{2}$. By H\"older's inequalities, we have to verify that:
$$\frac{1}{2}+\frac{1}{2q}-\frac{1}{N}-\frac{\e}{4q}-\frac{\e}{2N}+\frac{3\gamma-q}{q}+\frac{3q-6\gamma}{2q^{2}}\leq 1.$$
It means that we must have:
\begin{equation}
11q^{2}-6q(2+3\gamma)+18\gamma\geq0.
\label{egal2}
\end{equation}
We have next:
$$
\begin{aligned}
&\int^{t}_{0}\|\sqrt{f(s)\rho}u\cdot\n v_{1}\|_{L^{2}}^{2}ds,\\
&\leq \int^{t}_{0}\|\rho^{\frac{3\gamma-q}{q}}u\|_{L^{\frac{q}{3\gamma-q}}}^{2}\|\sqrt{f(s)}\D v_{1}\|^{2(1-\frac{\e}{2})}_{L^{p}}f(s)^{\frac{\e}{2}}\|\n v_{1}\|^{2\frac{\e}{2}}_{L^{2}}\|\rho\|_{L^{q}}^{\frac{3q-6\gamma}{q}}
ds,\\
&\leq \int^{t}_{0} \|\rho^{\frac{3\gamma-q}{q}}u\|_{L^{\frac{q}{3\gamma-q}}}^{2}(\|\rho\|_{L^{q}}^{\frac{1}{2}}
\|\sqrt{\rho f(s)}\p_{t}v_{1}\|_{L^{2}}+\sqrt{f(s)}\|\rho u\cdot \n v_{1}\|_{L^{p}}ds\\
&+\sqrt{f(s)}\|\rho u\cdot\n(\D)^{-1}\n( P(\rho))\|_{L^{p}}
+\sqrt{f(s)}\|\rho g\|_{L^{p}}+\sqrt{f(s)}\|\rho \p_{t}v\|_{L^{p}})^{2(1-\frac{\e}{2})}\\
&\hspace{9cm}\times\|\n v_{1}\|^{2\frac{\e}{2}}_{L^{2}(L^{2})}f(s)^{\frac{\e}{2}}\|\rho\|_{L^{q}}^{\frac{3q-6\gamma}{q}}ds,\\
&\leq \int^{t}_{0} \frac{1}{\theta}\|\rho^{\frac{3\gamma-q}{q}}u\|_{L^{\frac{q}{3\gamma-q}}}^{\frac{4}{\e}}f(s)\|\n v_{1}\|^{2}_{L^{2}}\|\rho\|_{L^{q}}^{\frac{2(3q-6\gamma)}{\e q}}+\theta\big(
\|\rho\|_{L^{q}}^{\frac{1}{2}}
\|\sqrt{\rho f(s)}\p_{t}v_{1}\|^{2}_{L^{2}}\\
&+\|\rho\|^{\frac{1}{2}}_{L^{q}}\|\sqrt{f(s)\rho} u\cdot \n v_{1}\|^{2}_{L^{2}}+\|\sqrt{f(s)}\rho u\cdot\n(\D)^{-1}\n( P(\rho))\|^{2}_{L^{p}}
+\|\sqrt{f(s)}\rho g\|^{2}_{L^{p}}\\
&\hspace{10cm}+\|\sqrt{f(s)}\rho \p_{t}v\|_{L^{p}}^{2}\big),\\
\end{aligned}
$$
where we have $\theta$ which depends of $t$ for the bootstrap.
\subsubsection*{Estimate on the term $\int^{t}_{0}\int_{\T^{N}}(\rho g-\frac{1}{\nu}\rho\p_{s}v)f(s)\p_{s}v_{1}dxds$}
For the term $\int^{t}_{0}\int_{\T^{N}}(\rho g-\frac{1}{\nu}\rho\p_{s}v)f(s)\p_{s}v_{1}dxds$, we proceed exactly as in the previous section.
\subsection*{Control on the norm $\rho\in L^{\infty}$}
Next we come back to equation (\ref{g28}) to get a control in norm $L^{\infty}$ on the density, more precisely we have:
\begin{equation}
\begin{aligned}
\log(\rho(t,x))\leq&\log(\|\rho_{0}\|_{L^{\infty}})+C\|(\D)^{-1}{\rm div} m_{0}\|_{L^{\infty}}+C\|(\D)^{-1}{\rm div}(\rho u)\|_{L^{\infty}}\\
&+C \int^{t}_{0}\frac{1}{\sqrt{s}}\|[u_{j},R_{i}R_{j}](\rho u_{i})(s,\cdot)\|_{L^{\infty}}ds.
\end{aligned}
 \label{28}
\end{equation}
We recall then from the previous section that $\sqrt{f(s)}\D v_{1}\in L^{2}(L^{p})$ with $\frac{1}{p}=\frac{1}{2}+\frac{1}{2q}$, by Gagliardo-Nirenberg inequality we have:
$$(f(s))^{1-\frac{\e}{2}}\|\n v_{1}\|^{2}_{L^{q_{1}}}\leq \|\sqrt{f(s)}\D v_{1}\|^{2(1-\frac{\e}{2})}_{L^{p}}\|\n v_{1}\|^{2\frac{\e}{2}}_{L^{2}},$$
with $\frac{1}{q_{1}}=\frac{1}{2}+\frac{1}{2q}-\frac{1}{N}-\frac{\e}{4q}-\frac{\e}{2N}$ and $q>\frac{N}{2}$.
We have then
by using the results of R. Coifman et al in \cite{1}:
\begin{equation}
\begin{aligned}
 &(f(s))^{\frac{1}{2}-\frac{\e}{4}}\|[(v_{1})_{j},R_{i}R_{j}](\rho u_{i})(s,\cdot)\|_{W^{1,\alpha}}\leq \|\n v_{1}\|_{L^{2}}^{\frac{\e}{2}}\|\sqrt{f(s)} \D v_{1}\|_{L^{p}}^{1-\frac{\e}{2}}\|\rho^{\frac{3\gamma-q}{q}}u\|_{L^{\frac{q}{3\gamma-q}}}
\\
&\hspace{10cm}\times\|\rho\|_{L^{\frac{q^{2}}{2q-3\gamma}}}
^{\frac{2\gamma-3q}{q}},
\end{aligned}
\label{infinal}
\end{equation}
with $\frac{1}{\alpha}=\frac{2q-3\gamma}{q^{2}}+\frac{3\gamma-q}{q}+\frac{1}{6}+\frac{1}{2q}+\frac{\e}{3}<\frac{1}{3}$
let $7q^{2}-3q(5+3\gamma)+18\gamma>0$. After a small calculus we obtain $q>\frac{15}{14}+\frac{9}{14}\gamma+\frac{3}{14}\sqrt{9\gamma^{2}-26\gamma+23}$. From (\ref{infinal}), we have $(f(s))^{\frac{1}{2}-\frac{\e}{4}}\|[(v_{1})_{j},R_{i}R_{j}](\rho u_{i})(s,\cdot)\|_{W^{1,\alpha}}\in L^{2}_{s}$.\\ 
We obtain then by Sobolev embedding:
$$|\int^{t}_{0}\frac{1}{f(s)^{\frac{1}{2}-\frac{\e}{4}}}f(s)^{\frac{1}{2}-\frac{\e}{4}}\|[(v_{1})_{j},R_{i}R_{j}](\rho u_{i})(s,\cdot)\|_{L^{\infty}}ds|\leq C \|\rho^{\frac{3\gamma-q}{q}}u\|_{L^{\infty}(L^{\frac{q}{3\gamma-q}})}\|\rho\|_{L^{\infty}(L^{\frac{q^{2}}{2q-3\gamma}})}
^{\frac{2\gamma-3q}{q}}.$$
We proceed similarly for the term $\|[(\Lambda^{-1}(P(\rho))_{j},R_{i}R_{j}](\rho u_{i})(s,\cdot)\|_{L^{1}_{s}(L^{\infty})}$. This term is crucial because it gives the value of $q$ that we must choose. Indeed we have the results of R. Coifman et al in \cite{1} if $q<3\gamma$:
$$\|[(\Lambda^{-1}(P(\rho))_{j},R_{i}R_{j}](\rho u_{i})(s,\cdot)\|_{W^{1,\beta}}\leq\|P(\rho)\|_{L^{\frac{\gamma}{q}}}\|\|\rho^{\frac{3\gamma-q}{q}}u\|_{L^{\frac{q}{3\gamma-q}}}
\|\rho\|_{L^{\frac{q^{2}}{2q-3\gamma}}},$$
with $\frac{1}{\beta}=\frac{2q-3\gamma}{q^{2}}+\frac{3\gamma-q}{q}+\frac{\gamma}{q}<\frac{1}{3}$ which implies
$3q^{2}-q(9\gamma+5)+6\gamma>0$.
However we see by the previous inequality that we must have $q>3\gamma$.\\
As $u_{0}\in L^{\infty}$, we can show that $\rho^{\frac{1}{p}}u\in L^{\infty}(L^{p})$ with $p$ arbitrarly large because $\lambda=0$. We have then $\rho^{\frac{1}{p}}u\in L^{\infty}(L^{p})$, $\rho^{1-\frac{1}{p}}\in L^{\infty}(L^{\frac{q}{1-\frac{1}{p}}})$. We conclude that by H\"older's inequalities that:
$$\frac{q}{\gamma}+\frac{1}{q}<\frac{1}{3},$$
let $q>3(\gamma+1)$.\\
By this way we control
$\log\rho\in L^{\infty}$ by the fact that $q>3(\gamma+1)$ and that (\ref{egal1}) and (\ref{egal2}) are verified. It means that $\rho\in L^{\infty}$ and $\frac{1}{\rho}\in L^{\infty}$, from theorem \ref{theo1} we have seen that the inequality (\ref{1.21}) and (\ref{1.211}) are preserved during the time if $\rho\in L^{\infty}$.
We can now assume that $(\rho,u)$ verify (\ref{1.21}) and (\ref{1.211}). We want now prove the uniqueness of this solution.
\subsection{Uniqueness}
You want prove now the result of uniqueness. To do it we want use the result of P. Germain in \cite{PG}. In the sequel we will note $(\rho_{1},u_{1})$ the solution of the theorem \ref{theo1} which exits on the time interval $[0,T_{0}]$.
We have shown that our solution check $\rho\in L^{\infty}(L^{\infty})$. By theorem \ref{theo1}, we obtain that our solution verify the following inequalities:
\begin{equation}
 \begin{aligned}
&\sup_{0<t\leq+\infty}\int_{\T^{N}}[\frac{1}{2}\rho(t,x)|u(t,x)|^{2}+|P(\rho(t,x))|+\sigma(t)|\n u(t,x)|^{2}dx\\
&+\sup_{0<t\leq +\infty}\int_{\T^{N}}[\frac{1}{2}\rho(t,x)f(t)^{N}(\rho|\dot{u}(t,x)|^{2}+|\n \omega(t,x)|^{2})dx\\
&+\int^{+\infty}_{0}\int_{\T^{N}}[|\n u|^{2}+f(s)\rho|\dot{u}|^{2}+|\omega|^{2})+\sigma^{N}|\n \dot{u}|^{2}]dxdt\\
&\hspace{9cm}\leq C(C_{0}+C_{f})^{\theta},
 \end{aligned}
\label{h51.21}
\end{equation}
and we obtain moreover:
\begin{equation}
\begin{cases}
 \begin{aligned}
&\sqrt{\rho}\p_{t}u\in L^{2}_{t}(L^{2}(\T^{N})),\\
&\sqrt{t}{\cal P}u\in L^{2}_{T}(H^{2}(\T^{N})),\\
&\sqrt{t}G=\sqrt{t}[(\lambda+2\mu){\rm div}u-P(\rho)]\in L^{2}_{T}(H^{1}(\T^{N})),\\
&\sqrt{t}\n u\in L^{\infty}_{T}(L^{2}(\T^{N})),
 \end{aligned}
\end{cases}
\label{a513}
\end{equation}
Now by the result of P. Germain in \cite{PG}, we are able to prove that $(\rho,u)=(\rho_{1},u_{1})$ on $[0,T_{0}]$. To see this we have just to verify that $(\rho_{1},u_{1}$ verify the conditions of the theorem 2.2 of \cite{PG}.
For simplicity we prove only the result for $N=3$. We know that $\n\rho_{1}\in L^{\infty}(B^{\e}_{N,1})\hookrightarrow L_{T_{0}}^{\infty}(L^{N})$ and $\n u_{1}\in L^{1}_{T_{0}}(L^{\infty})$.
The main thing is to see that $\sqrt{t} \dot{u}_{1}\in L^{2}_{T_{0}}(L^{N})$.
\\
We recall by Gagliardo-Nirenberg inequalities that:
$$\sqrt{t}\|\dot{u}_{1}\|_{L^{3}}\leq(t^{\frac{1}{4}-\frac{\e}{2}}\|\dot{u}_{1}\|_{L^{2}})^{\frac{1}{2}}(t^{\frac{3}{2}-\frac{\e}{2}}
\|\n \dot{u}_{1}\|_{L^{2}})^{\frac{1}{2}}t^{\frac{\e}{2}}.$$
From the inequalities (\ref{b1.22}), we deduce that $\sqrt{t}\dot{u}_{1}\in L^{2}(L^{3})$.\\
We have proved that our solution is unique on $[0,T_{0}]$, we have to see now what happen for $t\geq T_{0}$.
When  $t\geq T_{0}$, we have in this case by using inequalities (\ref{h51.21}),
\begin{equation}
 \begin{aligned}
&\sup_{T_{0}\leq t\leq +\infty}\int_{\T^{N}}[\frac{1}{2}\rho(t,x)|u(t,x)|^{2}+|P(\rho(t,x))|+\sigma(t)|\n u(t,x)|^{2}dx\\
&+\sup_{T_{0}\leq t\leq +\infty}\int_{\T^{N}}[\frac{1}{2}\rho(t,x)f(T_{0})^{2}(\rho|\dot{u}(t,x)|^{2}+|\n \omega(t,x)|^{2})dx\\
&+\int^{+\infty}_{T_{0}}\int_{\T^{N}}[|\n u|^{2}+f(T_{0}))\rho|\dot{u}|^{2}+|\omega|^{2})+f(T_{0})^{2}|\n \dot{u}|^{2}]dxdt\\
&\hspace{9cm}\leq C(C_{0}+C_{f})^{\theta},
 \end{aligned}
\label{u51.21}
\end{equation}
We can check easily that $v_{1}\in L^{1}_{loc}([T_{0},+\infty),B^{1+\e}_{N,\infty})$ so by proceeding similarly as subsubsection \ref{section61}, we obtain that $\rho\in L^{\infty}(\R,B^{1+\e}_{N,\infty})$. We get then that $\n u\in L^{1}(\R, L^{\infty})$ and we conclude by using again the result of P. Germain in \cite{PG}.


\begin{thebibliography}{}

\bibitem{AP}
H. ABIDI and M. PAICU. \'Equation de Navier-Stokes avec densit\'e et viscosit\'e variables dans l'espace critique. \textit{Annales de l'institut Fourier}, 57 no. 3 (2007), p. 883-917.
\bibitem{5BC}
H. BAHOURI and J.-Y. CHEMIN, \'Equations d'ondes quasilin\'eaires et
estimation de Strichartz, \textit{Amer. J. Mathematics}, 121, (1999), 1337-1377.
\bibitem{BCD}
H. BAHOURI, J.-Y. CHEMIN and R. DANCHIN. Fourier analysis and nonlinear partial differential equations,
\textit{to appear in Springer}.
\bibitem{5BJM}
J.-M. BONY, Calcul symbolique et propagation des singularit\'es pour
les \'equations aux d\'eriv\'ees partielles non lin\'eaires, \textit{Annales
Scientifiques de l'\'ecole Normale Sup\'erieure}, 14, (1981),
209-246.
\bibitem{aBD}
D. BRESCH and B. DESJARDINS, On the construction of approximate solutions for the 2d viscous shallow
water model and for compressible Navier-Stokes models. \textit{J. Maths Pures et Appl}., (86), 4, (2006), 262-268.
\bibitem{5BD}
D. BRESCH and B. DESJARDINS, Existence of global weak solutions for
a 2D Viscous shallow water equations and convergence to the
quasi-geostrophic model. \textit{Comm. Math. Phys.}, 238 (1-2), 211-223,
2003.
\bibitem{5BD1}
D. BRESCH and B. DESJARDINS, Existence of global weak solutions to
the Navier-Stokes equations for viscous compressible and heat
conducting fluids, \textit{Journal de Math\'ematiques Pures et Appliqu\'es},
Volume 87, Issue 1, January 2007, Pages 57-90.
\bibitem{5BD2}
D. BRESCH and B. DESJARDINS, Some diffusive capillary models of
Kortweg type. \textit{C. R. Math. Acad. Sci. Paris}, Section M\'ecanique,
332(11), 881-886, 2004.
\bibitem{5BDL}
D. BRESCH, B. DESJARDINS and C.-K. LIN, On some compressible fluid
models: Korteweg, lubrification and shallow water systems. \textit{Comm.
Partial Differential Equations}, 28(3-4), 843-868, 2003.
\bibitem{CKN}
L. CAFFARELLI, R. KOHN and L. NIRENBERG, Partial regularity of suitable weak solutions of the Navier-Stokes equations.,
\textit{Comm. Pure Appl. Math.}, 35(6): 771-831, 1982.
\bibitem{CD}
F. CHARVE and R. DANCHIN, Global existence in critical spaces for compressible Navier-Stokes equations, preprint.
\bibitem{5CT}
J.-Y. CHEMIN, Th\'eor\`emes d'unicit\'e pour le syst\`eme de
Navier-Stokes tridimensionnel, \textit{J. d'Analyse Math.} 77, (1999), 27-50.
\bibitem{5CL}
J.-Y. CHEMIN and N. LERNER, Flot de champs de vecteurs non
lipschitziens et \'equations de Navier-Stokes, \textit{J. Differential
Equations}, 121, (1992), 314-328.
\bibitem{CMZ}
Q. CHEN, C. MIAO and Z. ZHANG,
Global well-posedness for the compressible Navier-Stokes equations with the highly oscillating initial velocity, \textit{arXiv:0907.4540v2}.
\bibitem{CMZ1}
Q. CHEN, C. MIAO and Z. ZHANG, Well-posedness in critical spaces for the compressible
Navier-Stokes equations with density dependent viscosities, \textit{arXiv} :0811.4215v1 November 2008.
\bibitem{CCK2}
Y. CHO and H. KIM, Unique solvability for the density-dependent Navier-Stokes equations.
\bibitem{CCK3}
H. J. CHOE and H. KIM, Strong solution of the Navier-Stokes
equations for isentropic compressible fluids, \textit{J. Differential
Equations}, 190, (2003), 504-523.
\bibitem{5CK2}
H. J. CHOE and H. KIM, Strong solution of the Navier-Stokes
equations for nonhomogeneous incompressible fluids, \textit{Math. Meth.
Appl. Sci}. 28, (2005), 1-28.
\bibitem{1}
R. COIFMAN, P.-L. LIONS, Y. MEYER and S. SEMMES, Compensated-compactness and Hardy spaces, \textit{J. Math. Pures. Appl.}, 72 (1993), p247-286.
\bibitem{5CR}
F. COQUEL, D. DIEHL, C. MERKLE and C. ROHDE, Sharp and diffuse
interface methods for phase transition problems in liquid-vapour
flows. \textit{Numerical Methods for Hyperbolic and Kinetic Problems}, 239-270, IRMA Lect. Math. Theor.  Phys, 7, Eur. Math. Soc,
Z$\ddot{\mbox{u}}$rich, 2005.
\bibitem{DFourier}
R. DANCHIN, Fourier analysis method for PDE's, \textit{Preprint}, Novembre 2005.
\bibitem{DL}
R. DANCHIN, Local Theory in critical Spaces for Compressible Viscous
and Heat-Conductive Gases, \textit{Communication in
Partial Differential Equations}, 26 (78),1183-1233, (2001).
\bibitem{DG}
R. DANCHIN, Global Existence in Critical Spaces for Flows of
Compressible Viscous and Heat-Conductive Gases, \textit{Arch.Rational
Mech.Anal}.160, (2001), 1-39.
\bibitem{DU}
R. DANCHIN, On the uniqueness in critical spaces for compressible Navier-Stokes equations. \textit{NoDEA Nonlinear Differentiel Equations
Appl}, 12(1):111-128, 2005.
\bibitem{DW}
R. DANCHIN, Well-Posedness in Critical Spaces for Barotropic Viscous Fluids with Truly Not Constant
Density, \textit{Communications in Partial Differential Equations},32:9,1373-1397.
\bibitem{CCK5}
B. DESJARDINS, Regularity of weak solutions of the compressible isentropic Navier-Stokes equations, \textit{Comm Partial Differential Equations} 22 (1997) 977-1008.
\bibitem{5F}
E. FEIREISL, Dynmamics of Viscous Compressible Fluids,\textit{Oxford Lecture
Series in Mathematics and its Applications}, 26.
\bibitem{5F1}
E. FEIREISL. Compressible Navier-Stokes equations with a
non-monotone pressure law. \textit{J. Differential Equations},
184(1), 97-108, 2002.
\bibitem{5F2}
E. FEIREISL. On the motion of a viscous, compressible, and heat
conducting equation. \textit{Indiana Univ. Math. J.},
53(6), 1705-1738, 2004.
\bibitem{5F3}
E. FEIREISL, A. NOVOTNY and H. PETZELTOVA. On the
existence of globally defined weak solutions to
the Navier-Stokes equations. \textit{J. Math. Fluid Mech.}, 3(4), 358-392, 2001.
\bibitem{PG}
P. GERMAIN, Weak stong uniqueness for the isentropic compressible Navier-Stokes system, \textit{preprint}.
\bibitem{aGJX}
Z. GUO, Q. JIU and Z. XIN. Spherically symmetric isentropic compressible flows with density-dependent viscosity coefficients. \textit{SIAM J. Math. Anal.} (2008), 1402-1427.
\bibitem{H}
B. HASPOT, Cauchy problem for viscous shallow-water equations with a term of capillarity, accepted in \textit{M3AS}.
\bibitem{H1}
B. HASPOT, Cauchy problem for viscous shallow water equations  with a term of capillarity , accepted in \textit{Hyperbolic Problems: Theory, Numerics and Application} PSAPM volume, edited by Drs. Tadmor, Liu, Tzavaras.
\bibitem{H2}
B. HASPOT, Local well-posedness results for density-dependent incompressible fluids, \textit{Arxiv},:0902.1982 (February 2009).
\bibitem{H3}
B. HASPOT, Well-posedness in critical spaces for barotropic viscous fluids , \textit{Arxiv},:0902.1982 (February 2009).
\bibitem{Hoffuni}
D. HOFF. Uniqueness of weak solutions of the Navier-Stokes equations of multidimensional compressible flows.
\textit{SIAM. J. Anal.},vol 37, No6, 1742-1760 (2006).
\bibitem{Hoffn1}
D. HOFF. Compressible flow in a Half-space with Navier boundary conditions, \textit{Journal of Mathematical Fluid Mechanics}, Volume 7, Number 3 / August, (2005) 315-338.
\bibitem{Hoffn2}
D. HOFF. Dynamics of singularity surfaces for compressible, viscous flows in two space dimensions, \textit{Communications on Pure and Applied Mathematics},
Volume 55 Issue 11, 1365 - 1407 (2002).
\bibitem{9}
D. HOFF. Global existence for 1D, compressible, isentropic
Navier-Stokes equations with large initial data.
\textit{Trans. Amer. Math. Soc}, 303(1), 169-181, 1987.
\bibitem{10}
D. HOFF. Global well-posedness of the Cauchy problem for nonisentropic gas dynamics with discontinuous data, \textit{J. Diff. Eq.}, 95 (1992), p33-73
\bibitem{5H4}
D. HOFF. Discontinuous solutions of the Navier-Stokes equations
for multidimensional flows of the heat conducting fluids.
\textit{Arch. Rational Mech. Anal.}, 139, (1997), p. 303-354.
\bibitem{5H2}
D. HOFF. Global solutions of the Navier-Stokes equations for
multidimensional compressible flow with discontinuous initial
data. \textit{J. Differential Equations}, 120(1), 215-254, 1995.
\bibitem{5H3}
D. HOFF. Strong convergence to global solutions for
multidimensional flows of compressible, viscous fluids with
polytropic equations of state and discontinuous initial data.
\textit{Arch. Rational Mech. Anal.}, 132(1), 1-14, 1995.
\bibitem{5HZ}
D. HOFF and K. ZUMBRUM. Multi-dimensional diffusion waves for
the Navier-Stokes equations of compressible flow,
\textit{Indiana University Mathematics Journal}, 1995, 44, 603-676.
\bibitem{5J}
S. JIANG and P. ZHANG. Axisymetrics solutions of the 3D
Navier-Stokes equations for compressible isentropic fluids.
J. Math. Pures Appl. (9), 82(8): 949-973, 2003.
\bibitem{5K1}
A. V. KAZHIKOV. The equation of potential flows of a compressible
viscous fluid for small Reynolds numbers: existence,
uniqueness and stabilization of solutions. \textit{Sibirsk. Mat. Zh.}, 34 (1993), no. 3, p. 70-80.
\bibitem{11}
A. V. KAZHIKOV and V. V. SHELUKHIN. Unique global solution with
respect to time of initial-boundary value problems for one-
dimensional equations of a viscous gas. \textit{Prikl. Mat. Meh.}, 41(2): 282-291, 1977.
\bibitem{aJJX}
H.-L. LI, J. LI and Z. XIN. Vanishing of vacuum states and blowup phenomena of the compressible Navier-Stokes equation. \textit{Accepted in Comm. Math. Phys}, (2007).
\bibitem{13}
P.-L. LIONS. Mathematical Topics in Fluid Mechanics, Vol 2,
Compressible models, \textit{Oxford University Press}, (1998).
\bibitem{5MN}
A. MATSUMURA and T. NISHIDA. The initial value problem for
the equations of motion of compressible
viscous and heat-conductive fluids. \textit{Proc. Japan Acad. Ser. A Math. Sci}, 55(9), 337-342, 1979.
\bibitem{5MV}
A. MELLET and A. VASSEUR, On the isentropic compressible Navier-Stokes
equation, Arxiv preprint math.AP/0511210, 2005 - arxiv.org
\bibitem{5MV2}
A. MELLET and A. VASSEUR, $L^{p}$ estimates for quantities advected by a compressible flow, \textit{Journal of Mathematical Analysis and Applications}, Volume 355, Issue 2, 548-563 (2009).
\bibitem{aMV1}
A. MELLET and A. VASSEUR, Existence and uniqueness of global strong solutions for one-dimensional compressible Navier-Stokes equations. \textit{SIAM J. Math. Anal.} 39 (2007/2008), no. 4, 1344-1365.
\bibitem{5Na}
J. NASH, Le probl\`eme de Cauchy pour les \'equations
diff\'erentielles d'un fluide g\'en\'eral, \textit{Bulletin de
la Soci\'et\'e Math\'ematique de France}, 1962, 90, 487-497.
\bibitem{NS}
A. NOVOTNY\'{Y} and I. STRA\v{S}KRABA, Introduction to the mathematical theory of compressible flow, \textit{Oxford lecture series in mathematics and its applications}. \textbf{27}.
\bibitem{5Ro}
C. ROHDE, On local and non-local Navier-Stokes-Korteweg systems for
liquid-vapour phase transitions. \textit{ZAMMZ. Angew. Math. Mech.}, 85(2005), no. 12, 839-857.
\bibitem{5RS}
T. RUNST and W. SICKEL. Sobolev spaces of fractional order, Nemytskij operators, and nonlinear partial differential equations. de Gruyter Series in
Nonlinear Analysis and Applications, 3. \textit{Walter de Gruyter and Co.}, Berlin (1996).
\bibitem{CCK17}
R. SALVI and I. STRASKRABA, Global existence for viscous compressible fluids and their behavior at $t\rightarrow+\infty$, \textit{J. Fac. Sci. Univ. Tokyo Sect.}IA, Math, 40 (1993) 17-51.
\bibitem{16}
D. SERRE. Solutions faibles globales des \'equations de
Navier-Stokes pour un fluide compressible.,
303(13): 639-642, 1986
\bibitem{5So}
V. A. SOLONNIKOV. Estimates for solutions of nonstationary
Navier-Stokes systems. \textit{Zap. Nauchn. Sem. LOMI}, 38,
(1973), p.153-231; J. Soviet Math. 8, (1977), p. 467-529.
\bibitem{19}
V. A. SOLONNIKOV. Solvability of the initial boundary value problem for the equation of a viscous compressible fluid, \textit{J. Sov. Math.}, 14 (1980), p1120-1133.
\bibitem{20}
V. VALLI, W. ZAJACZKOWSKI. Navier-Stokes equations for compressible
fluids: global existence and qualitative properties of the solutions
in the general case. \textit{Commun. Math. Phys.}, 103, (1986) no 2, p.
259-296.
\bibitem{aVD}
J. D. VAN DER WAALS. Thermodynamische Theorie der Kapillarit\"at unter Voraussetzung stetiger
Dichte\u{a}nderung, \textit{Z. Phys. Chem.} 13, (1894), 657-725.
\bibitem{V6}
A. VASSEUR. A new proof of partial regularity of solutions to Navier-Stokes equations. \textit{NoDEA Nonlinear Differential Equations Appl.}, 14(5-6):753-785, 2007.
\bibitem{CCK25}
V. A. WEIGANT, Example of non-existence in the large for the problem of the existence of solutions of Navier-Stokes equations for compressible viscous barotropic fluids, \textit{Dokl. Akad. Nauk} 339 (1994) 155-156 (in Russian).
\bibitem{5W}
W. WANG and C-J XU. The Cauchy problem for viscous shallow
water equations.
\textit{Rev. Mat. Iberoamericana}  21., no. 1 (2005), 1-24.
\end{thebibliography}
\end{document}